\documentclass{amsart}
\usepackage[all]{xy}
\usepackage{amsmath}
\usepackage{amssymb}
\usepackage{amscd}
\usepackage{amsthm}
\usepackage{graphicx}
\usepackage[T1]{fontenc}
\usepackage{multirow}
\usepackage[francais]{babel}

\newcommand{\F}{\mathbb{F}}
\newcommand{\Fq}{\mathbb{F}_q}
\newcommand{\Q}{\mathbb{Q}}
\newcommand{\barQl}{{\bar\Q_\ell}}

\newcommand{\Z}{\mathbb{Z}}

\newcommand{\cE}{\mathcal{E}}
\newcommand{\cF}{\mathcal{F}}
\newcommand{\cH}{\mathcal{H}}

\newcommand{\cL}{\mathcal{L}}
\newcommand{\cM}{\mathcal{M}}
\newcommand{\cN}{\mathcal{N}}
\newcommand{\cO}{\mathcal{O}}

\newcommand{\bG}{\mathbf{G}}
\newcommand{\bH}{\mathbf{H}}
\newcommand{\bK}{\mathbf{K}}
\newcommand{\bL}{\mathbf{L}}
\newcommand{\bM}{\mathbf{M}}
\newcommand{\bP}{\mathbf{P}}
\newcommand{\bS}{\mathbf{S}}
\newcommand{\bT}{\mathbf{T}}

\newcommand{\GL}{\mathrm{GL}}
\newcommand{\SL}{\mathrm{SL}}
\newcommand{\Spin}{\mathrm{Spin}}
\newcommand{\Or}{\mathrm{O}}
\newcommand{\SO}{\mathrm{SO}}
\newcommand{\hSpin}{\text{$\frac{1}{2}$}\Spin}
\newcommand{\PSO}{\mathrm{PSO}}
\newcommand{\Sp}{\mathrm{Sp}}
\newcommand{\PSp}{\mathrm{PSp}}
\renewcommand{\AA}{{\mathrm{A}}}
\newcommand{\BB}{{\mathrm{B}}}
\newcommand{\CC}{{\mathrm{C}}}
\newcommand{\DD}{{\mathrm{D}}}
\newcommand{\EE}{{\mathrm{E}}}
\newcommand{\FF}{{\mathrm{F}}}
\newcommand{\GG}{{\mathrm{G}}}

\newcommand{\ad}{{\mathrm{ad}}}
\newcommand{\der}{{\mathrm{der}}}
\newcommand{\ext}{{\mathrm{ex}}}
\newcommand{\simc}{{\mathrm{sc}}}
\newcommand{\semi}{{\mathrm{ss}}}
\newcommand{\cusp}{{0}}
\newcommand{\triv}{k}

\DeclareMathOperator{\ind}{ind}

\DeclareMathOperator{\Irr}{Irr}
\DeclareMathOperator{\IC}{IC}
\DeclareMathOperator{\End}{End}
\DeclareMathOperator{\Hom}{Hom}

\DeclareMathOperator{\Tr}{Tr}
\DeclareMathOperator{\Norm}{N}
\newcommand{\Cent}{\mathrm{Z}}
\newcommand{\Cento}{\Cent^\circ}

\DeclareMathOperator{\pgcd}{pgcd}

\def\isom{\overset\sim\to}
\def\cf{{{\it cf.\ }}}
\def\ie{{\it i.e.,\ }}

\newcommand{\tE}{\widetilde E}

\newcommand{\unu}{{\underline \nu}}
\newcommand{\uw}{{\underline w}}
\newcommand{\WGSE}{{W^\bG_{\Sigma,\cE}}}
\newcommand{\tWGSE}{{\widetilde W^\bG_{\Sigma,\cE}}}
\newcommand{\WGL}{{W^\bG_\bL}}
\newcommand{\WHM}[1][{}]{{W^\bH_{\bM_{#1}}}}
\newcommand{\WHMj}{\WHM[j]}
\newcommand{\tWHMj}{{\widetilde W^\bH_{\bM_j}}}

\newcommand{\Gsigma}{{\bH}}

\newtheorem{thm}{Th\'eor\`eme}[section]

\newtheorem{lem}[thm]{Lemme}
\newtheorem{cor}[thm]{Corollaire}
\theoremstyle{definition}

\newtheorem{remark}[thm]{Remarque}
\theoremstyle{remark}
\newtheorem{exm}[thm]{Exemple}

\catcode`\@=11
\def\@@and{et}
 \def\courriel#1#2{\begingroup
     \@ifnotempty{#2}{\nobreak\indent{\itshape Adresse courriel}%
       \@ifnotempty{#1}{, \ignorespaces#1\unskip} :\space
       \ttfamily#2\par}\endgroup}
 \renewcommand{\email}[2][]{\g@addto@macro\addresses{\courriel{#1}{#2}}}
\catcode`\@=12

\title{Localisation de faisceaux caract\`eres}
\author{Pramod N.~Achar}
\thanks{Le premier auteur \'etait partiellement appuy\'e par la subvention
DMS-0500873 de la NSF}
\address{Department of Mathematics\\
  Louisiana State University\\
  Baton Rouge, LA 70803, USA}
\email{pramod@math.lsu.edu}
\author{Anne-Marie Aubert}
\address{C.N.R.S.\\Institut de Math\'ematiques de Jussieu\\
175 rue du Chevaleret, 75013 Paris
  \\
  France}
\email{aubert@math.jussieu.fr}
\date{\today}

\keywords{faisceaux caract\`eres, caract\`eres des groupes
r\'eductifs finis, repr\'esentations des groupes de Weyl, classes
unipotentes, support unipotent}

\subjclass{20C33, 20G40}
\begin{document}

\begin{abstract}
Nous obtenons une formule pour les valeurs de la fonction
caract\'eristique d'un faisceau caract\`ere en fonction de la th\'eorie des
repr\'esentations de certains groupes
finis, li\'es au groupe de Weyl.  Cette formule, qui g\'en\'eralise des
r\'esultats ant\'erieurs de M\oe glin et de Waldspurger, d\'epend
de la connaissance de certains sous-groupes r\'eductifs admettant un
faisceau caract\`ere cuspidal.   Dans un
second temps, afin de rendre la formule plus
explicite dans le cas d'un groupe quasi-simple, nous d\'eterminons ces
sous-groupes \`a conjugaison pr\`es.

\bigskip

\noindent
\textsc{Abstract.}  We obtain a
formula for the values of the characteristic function of a character
sheaf, in terms of the representation theory of certain finite groups
related to the Weyl group.  This formula, a generalization of previous
results due to M\oe glin and Waldspurger, depends on knowledge
of certain reductive subgroups that admit cuspidal character sheaves.  For
quasi-simple groups, we make the formula truly explicit by determining all
such subgroups upto conjugation.
\end{abstract}

\maketitle

\section{Introduction}
\label{sect:intro}
Soit $\bG$ un groupe alg\'ebrique r\'eductif connexe sur la cl\^oture
alg\'ebrique d'un corps fini $\Fq$ et d\'efini sur ce dernier. Nous noterons
$G=\bG^F$ le groupe (fini) des \'el\'ements de $\bG$ fix\'es par 
l'endomorphisme de Frobenius $F$ associ\'e \`a la
structure $\F_q$-rationnelle du groupe $\bG$. 

Soit $\bL$ un sous-groupe de L\'evi d'un sous-groupe parabolique de $\bG$.  
Nous notons $\bL_{\der}$ et $\Cento(\bL)$ respectivement son groupe d\'eriv\'e 
et la composante neutre de son centre.
Nous consid\'erons une classe de conjugaison
$\cO$ de $\bL/\Cento(\bL)$, nous notons $\Sigma$ l'image r\'eciproque de celle-ci
dans $\bL$ via la projection naturelle $\bL\to \bL/\Cento(\bL)$
et $\Sigma_\semi$ l'ensemble des parties semi-simples des
\'el\'ements de $\Sigma$. Nous supposons l'existence d'un syst\`eme local
$\cE$ sur $\Sigma$ tel que la paire $(\Sigma,\cE)$ soit cuspidale au sens
de \cite[2.4]{LIC}.

Soit $g$ un \'el\'ement de $G$, fix\'e une fois pour toutes. Nous
l'\'ecrivons $g=\sigma
v$, o\`u $\sigma$ est semi-simple et $v$ est unipotent et commute \`a
$\sigma$.
Nous posons $\bH:=\Cento_\bG(\sigma)$.

Soit $x$ un \'el\'ement de $\bG$ tel que $x^{-1}\sigma x\in \Sigma_\semi$. On a 
alors $\sigma\in x \bL x^{-1}$.
Posons
\[ \bM_x:=(x\bL x^{-1})\cap\bH=\Cent_\bH(x\Cento(\bL) x^{-1})  \]
et consid\'erons l'ensemble
\[ \cM = \{ \bM_x \,:\, x^{-1} \sigma x \in \Sigma_\semi \} \]
de sous-groupes de L\'evi de $\bH$ (isomorphes car conjugu\'es sous
$\bG$).
Nous dirons que deux \'el\'ements de $\cM$ sont \'equivalents s'ils sont
conjugu\'es sous $\bH$.
L'ensemble $\cM$ est ainsi partitionn\'e en un nombre fini de classes 
d'\'equivalence \[ \cM = \cM_1 \sqcup \cdots \sqcup \cM_r. \]
Nous d\'emontrerons au Th\'eor\`eme~\ref{thm:m-conj}
que, si $\bG$ est semi-simple, quasi-simple, et diff\'erent de $\PSp_{2n}$,
$\PSO_{2n}$, $\hSpin_{2n}$ et $\EE_7^\simc$, alors tous les membres de
$\cM$ sont conjugu\'es sous $\bH$ et que si $\bG$ est l'un de ces quatres
groupes, alors $\cM$ se r\'epartit en une ou deux classes de conjugaison
sous $\bH$.

Pour tout $j \in \{1, \ldots, r\}$, nous fixons un \'el\'ement $g_{a_j} \in \bG$ 
tel que $\bM_{g_{a_j}} \in \cM_j$, et nous posons
\[
\bT^\bH_j := g_{a_j} \Cento(\bL) g_{a_j}^{-1}
\qquad\text{et}\qquad
\bM_j := \bM_{g_{a_j}} = \Cent_\bH(\bT^\bH_j).
\]
Nous posons
\[\WHMj :=\Norm_{\bH}(\bM_j)/\bM_j\quad\text{et}\quad \WGL:=\Norm_\bG(\bL)/\bL.\]
Soit $a_j$ l'\'el\'ement de $\WGL$ tel que $F(a_j)$ 
soit l'image de $g_{a_j}^{-1} F(g_{a_j})$ dans $\WGL$.
L'application $\tilde\iota_j\colon\Norm_\bH(\bM_j)\to\Norm_\bG(\bL)$
d\'efinie  par $\tilde\iota_j(h)=g_{a_j}^{-1}hg_{a_j}$
induit un plongement
\[
\iota_j\colon \WHMj\hookrightarrow \WGL.
\]
Le morphisme de Frobenius agit \`a la fois sur $\WGL$ et sur chacun des 
$\WHMj$, mais
le plongement $\iota_j$ n'est 
en g\'en\'eral pas \'equivariant pour ces deux actions. Nous garderons la 
notation $F\colon \WGL \to \WGL$ pour l'automorphisme de $\WGL$ induit par le
Frobenius, et noterons $\eta_j\colon\WHMj \to \WHMj$ l'automorphisme
\emph{inverse} \`a celui qui est induit par le Frobenius sur $\WHMj$.

Soit $\tWHMj$ le produit semi-direct de $\WHMj$ par le groupe cyclique
engendr\'e par $\eta_j$. Nous notons $\Irr(\WHMj)_\ext$ l'ensemble des
repr\'esentations irr\'eductibles de $\WHMj$ qui s'\'etendent en des 
repr\'esentations de $\tWHMj$.
Pour tout $E' \in \Irr(\WHMj)_\ext$, nous fixons une fois pour toutes une
extension $\tE'$ \`a $\tWHMj$.

Nous associons \`a tout \'el\'ement $w$ de $\WHMj$ une certaine fonction
de Green $Q_w$ (voir~(\ref{eqn:Qw})) et pour tout $E' \in \Irr(\WHMj)_\ext$, nous
posons
\[
Q_{E'}(v) := \frac{1}{|\WHMj|} \sum_{w' \in \WHMj} \Tr(\eta_j w', \tilde
E') Q_{w'}(v).
\]

\smallskip

Nous prenons pour $\cE$ (syst\`eme local sur $\Sigma$) l'image r\'eciproque,
sous l'application naturelle $\bL\to \bL/\bL_{\der}\times \bL/\Cento(\bL)$, 
du produit tensoriel externe d'un syst\`eme local Kummerien de rang
$1$ sur le tore $\bL/\bL_{\der}$ et d'un syst\`eme local irr\'eductible
$\bL$-\'equivariant sur $\cO$. 
Nous supposons que la paire cuspidale $(\Sigma,\cE)$ est $F$-stable, et 
fixons un isomorphisme $\varphi_0\colon F^*\cE\isom\cE$.
Nous posons
\[\WGSE:=\left\{n\in\Norm_\bG(\bL)\,:\,n\Sigma
n^{-1}=\Sigma,\;\ad(n)^*\cE\isom\cE\right\}/\bL.\]
Nous supposons que $\WGSE$ est produit semi-direct d'un groupe de
Coxeter fini par un groupe ab\'elien fini (dans le cas o\`u le centre de
$\bG$ est connexe, $\WGSE$ est un groupe de Coxeter fini, voir par exemple
\cite[(5.16)]{ShFC}, \cite[\S~4.2]{ShFCII}).
Nous notons $\Phi^+_{\Sigma,\cE}$
l'ensemble des racines positives associ\'ees \`a ce groupe de Coxeter.
Soit \[Z_{\Sigma,\cE}=\left\{n\in\Norm_\bG(\bL)\,:\,F(n\Sigma n^{-1})=\Sigma,
\;\ad(n)^*F^*\cE\isom\cE\right\}/\bL.\]
Il existe un \'el\'ement unique $w_1$ du groupe de Weyl de $\bG$ tel que
$Z_{\Sigma,\cE}=w_1 W_{\Sigma,\cE}^\bG$ et tel que l'application $\gamma_1\colon 
\WGSE \to \WGSE$ d\'efinie par $\gamma_1(w) =
w_1^{-1}F^{-1}(w) w_1$ envoie tout \'el\'ement de
$\Phi_{\Sigma,\cE}^+$
sur une racine positive.
Soit $\tWGSE$ le produit semi-direct de $\WGSE$
par le groupe cyclique engendr\'e par $\gamma_1$. Nous notons
$\Irr(W_{\Sigma,\cE}^\bG)_{\ext}$ l'ensemble des repr\'esentations
irr\'eductibles
de $W_{\Sigma,\cE}^\bG$ qui s'\'etendent en des repr\'esentations de
$\tWGSE$. Pour $E\in\Irr(W_{\Sigma,\cE}^\bG)_{\ext}$, nous choisissons
une extension $\tE$ de $E$ qui est d\'efinie sur $\Q$.

\smallskip

Un r\^ole important sera jou\'e par les ensembles de double classes
\[ \cN_j := \WGSE \backslash \WGL / \WHMj,  \]
pour $1\le j\le r$.
Dans chaque double classe $\unu \in \cN_j$, nous choisissons, une fois pour 
toutes, un repr\'esentant $w_\unu$, nous posons 
\[
W(\unu) = w_1 \WGSE \cap F^{-1}(w_\unu) (a_j \WHM_j) w_\unu^{-1},
\]
et nous d\'efinissons deux plongements comme suit:
\[
\begin{aligned}
\lambda\colon& W(\unu) \to \WHMj, \qquad& \lambda(\uw) &= \eta_j^{-1}(a_j^{-1}
F^{-1}(w_\unu^{-1}) \uw w_\unu) \\
\kappa\colon& W(\unu) \to \WGSE, \qquad& \kappa(\uw) &= \gamma_1^{-1}(w_1^{-1}
\uw).
\end{aligned}
\]
Nous introduisons maintenant une famille d'accouplements, param\'etr\'ee par les $\cN_j$, qui relient les ensembles
$\Irr(\WGSE)_\ext$ et $\Irr(\WHMj)_\ext$.  Soit
$\unu \in
\cN_j$.  Si $E \in \Irr(\WGSE)_\ext$ et $E' \in
\Irr(\WHMj)_\ext$, on pose
\[
\langle E, E' \rangle_\unu 
:= \frac{1}{|W(\unu)|} \sum_{\uw \in W(\unu)} \Tr(\gamma_1 \kappa(\uw), \tilde E)
\Tr(\eta_j \lambda(\uw), \tilde E').
\]
Nous choisissons un
repr\'esentant $\dot w_1$ de $w_1$ dans $\Norm_\bG(\bL)$ ainsi qu'un 
\'el\'ement $g_{\dot w_1}$
de $\bG$ tel que $g_{\dot w_1}^{-1}F(g_{\dot w_1})=F(\dot w_1)$, et nous posons
\[\bL^{w_1}=g_{\dot w_1}\bL g_{\dot w_1}^{-1},\quad
\Sigma^{w_1}=g_{\dot w_1}\Sigma g_{\dot w_1}^{-1},\quad \cE^{w_1}=\ad(g_{\dot
w_1}^{-1})^*\cE.
\]
Nous posons \[A_0=\IC(\bar\Sigma,\cE)[\dim\Sigma]\] (lequel est un faisceau
caract\`ere cuspidal sur $\bL$, si la caract\'eristique de $\Fq$ est bonne pour
$\bG$) et nous notons 
\[K=K(\bL,\Sigma,\cE)=\ind_\bP^\bG A_0\] le faisceau pervers semi-simple sur $\bG$
paraboliquement induit \`a partir de $A_0$ au sens de Lusztig; $K$
est naturellement muni d'une
structure mixte, nous notons $\varphi\colon
F^*K\isom K$ l'isomorphisme correspondant \`a cette derni\`ere et
$\chi_{K,\varphi}\colon G\to\barQl$ la fonction caract\'eristique de
$(K,\varphi)$. 
Les composantes irr\'eductibles de $K$ sont des faisceaux caract\`eres
$F$-stables (la structure mixte $\varphi_A\colon F^*A\isom A$ sur un tel
faisceau
caract\`ere $A$ est induite par $\varphi$) et tout faisceau caract\`ere
$F$-stable sur $\bG$ est une composante d'un induit de ce type.

Nous d\'efinissons alors \[K^{w_1}:=K(\bL^{w_1},\Sigma^{w_1},\cE^{w_1}),\] 
nous fixons
$\varphi_0^{w_1}\colon F^*\cE^{w_1}\isom\cE^{w_1}$ de sorte que
l'isomorphisme $F^* A_0\isom A_0$ induit co\"\i ncide avec
$\varphi_{A_0}$ et nous notons
$\varphi^{w_1}\colon F^* K^{w_1}\isom K^{w_1}$ l'isomorphisme 
par $\varphi$.

Le complexe $K^{w_1}$ admet la d\'ecomposition suivante:
\[K^{w_1}=\bigoplus_A A\otimes V_A,\]
o\`u $V_A=\Hom_{\cM G}(A,K^{w_1})$ est une repr\'esentation irr\'eductible
de $W_{\Sigma,\cE}^\bG$.

Pour tout faisceau caract\`ere $A=A_E$, avec $E \in \Irr(\WGSE)_\ext$,
nous obtenons au Th\'eor\`eme~\ref{thm:princ}
la formule suivante pour
la valeur de la fonction caract\'eristique
de $A_E$ \[
\chi_{A_E}(\sigma v) = \sum_{j=1}^r \sum_{\unu \in \cN_j}\, \sum_{E' \in
\Irr(\WHMj)_\ext} \langle E,E'\rangle_\unu\, Q_{E'}(v).
\]
Un cas particulier de la formule ci-dessus (celui correspondant \`a $\sigma=1$),
d\^u \`a Lusztig (\cf \cite{LCV}) a \'et\'e l'un des ingr\'edients
essentiels de \cite{AA}. Dans le cas ``oppos\'e'' au pr\'ec\'edent
(correspondant aux faisceaux caract\`eres dans la s\'erie unipotente), une
formule du type ci-dessus a \'et\'e obtenue par Shoji en
\cite[Lemma~4.5]{Sh}. Des formules g\'en\'erales pour les groupes symplectiques et
sp\'eciaux orthogonaux figurent dans les travaux de M\oe glin et
Waldspurger (l'entier $r$ est alors \'egal \`a $1$). Notre formule en est 
inspir\'ee (en particulier de
\cite[Proposition~7.2]{W}) et notre d\'emonstration
est une combinaison de \cite[Proposition~2.16]{MW} et \cite[Lemme~7.1]{W}. 
Notre formule est cependant moins explicite que dans {\it loc. cit.} dans
la mesure o\`u nous n'avons pas explicit\'e les structures mixtes
concern\'ees.

\section{Quelques Rappels sur les Faisceaux Caract\`eres}
\label{sect:rappels}

\subsection{Complexes admissibles}

Soient $\bar\F_q$ la cl\^oture alg\'ebrique d'un corps fini $\F_q$ de
caract\'eristique not\'ee $p$ et
$\bG$ un groupe alg\'ebrique r\'eductif connexe sur $\bar\F_q$ qui est d\'efini 
sur $\F_q$. Nous noterons $F$ l'endomorphisme de Frobenius associ\'e \`a la 
structure $\F_q$-rationnelle de $\bG$ et $G$ le groupe (fini) $\bG^F$ des points de
$\bG$ fixes par $F$. 

Nous notons $\cM(\bG)$ la cat\'egorie des faisceaux pervers sur $\bG$.

Soit $\bL$ un sous-groupe de L\'evi d'un sous-groupe parabolique $\bP$ de
$\bG$.
Nous notons $\bL_{\der}$ le groupe d\'eriv\'e de $\bL$ et $\bT_1$ 
la composante neutre $\Cento(\bL)$ du centre de $\bL$
(donc $\bL=\Cent_\bG(\bT_1)$). 

Soit $\Sigma$ l'image r\'eciproque dans $\bL$ d'une classe de conjugaison
$\cO$ de $\bL/\bT_1$ sous la projection naturelle
$\bL\to \bL/\bT_1$.
Soit $\cE$ un syst\`eme local sur $\Sigma$, qui est l'image r\'eciproque,
sous l'application naturelle $\bL\to \bL/\bL_{\der}\times \bL/\bT_1$, de
$\cL'\boxtimes\cE'$, o\`u $\cL'$ est un syst\`eme local Kummerien de rang
$1$ sur le tore $\bL/\bL_{\der}$ et $\cE'$ est un syst\`eme local irr\'eductible
$\bL$-\'equivariant (pour l'action de conjugaison) sur $\cO$. 

Nous supposons que la paire $(\Sigma,\cE)$ est cuspidale au sens de
\cite[2.4]{LIC}. Nous la supposons aussi $F$-stable (\ie
$F(\Sigma)=\Sigma$ et
$F^*\cE\isom\cE$), nous fixons un isomorphisme $\varphi'\colon 
F^*\cE'\isom \cE'$ et notons $\varphi_0\colon
F^*\cE\isom\cE$ l'isomorphisme induit par $\varphi'$.
Nous posons $A_0=\IC(\bar\Sigma,\cE)[\dim\Sigma]$. 

\begin{remark} Si la caract\'eristique $p$ est
\emph{presque bonne} pour $\bG$ (\ie $p$ est bonne pour tout facteur de $\bG$
de type exceptionnel et il n'y a pas de condition pour les facteurs de
type classique), alors $A_0$ est un faisceau caract\`ere cuspidal sur $\bL$.
\end{remark}

\smallskip
Soit $K=K(\bL,\Sigma,\cE)=\ind_\bP^\bG A_0$ le faisceau pervers sur $\bG$
induit \`a partir de $A_0$ (\cf \cite[\S 4.1]{LCS}). Il est semi-simple et est
naturellement muni d'une
structure mixte (\cf \cite[\S 8.1]{LCS}). Nous notons $\varphi\colon
F^*K\isom K$ l'isomorphisme correspondant et
$\chi_{K,\varphi}\colon G\to\barQl$ la fonction caract\'eristique de
$(K,\varphi)$, qui est une fonction centrale sur $G$ et est d\'efinie
par
\[\chi_{K,\varphi}(x)=\sum_i(-1)^i\Tr(\varphi,\cH^i_x(K)),\]
o\`u $\cH^i_x(K)$ d\'esigne la fibre en $x\in G$ du $i$-\`eme faisceau de
cohomologie $\cH^i(K)$ de $K$.

Les composantes irr\'eductibles de $K$ sont des faisceaux caract\`eres
$F$-stables (la structure mixte $\varphi_A\colon F^*A\isom A$ sur un tel
faisceau
caract\`ere $A$ est induite par $\varphi$) et tout faisceau caract\`ere
$F$-stable sur $\bG$ est composante d'un induit de ce type. 

Nous posons
\begin{equation} \label{WGSE}
\WGSE:=\left\{n\in\Norm_\bG(\bL)\,:\,n\Sigma
n^{-1}=\Sigma,\;\ad(n)^*\cE\isom\cE\right\}/\bL.
\end{equation}

L'alg\`ebre d'endomorphismes $\End_{\cM \bG}(K)$ de $K$ dans $\bG$ est
isomorphe \`a l'alg\`ebre de groupe $\barQl \WGSE$ tordue par
un
$2$-cocyle (\cite[3.4]{LIC}). Shoji a montr\'e en \cite[Lem. 5.9]{ShFC}
que le
cocycle est trivial lorsque le centre $\Cent(\bG)$ de $\bG$ est connexe et
que le groupe
$\bG/\Cent(\bG)$ est
simple. \emph{Nous supposons dor\'enavant le cocycle trivial}.

\smallskip
Soit
\[
Z_{\Sigma,\cE}=\left\{n\in\Norm_\bG(\bL)\,:\,F(n\Sigma n^{-1})=\Sigma,
\;\ad(n)^*F^*\cE\isom\cE\right\}/\bL.
\]
Lorsque le centre de $\bG$ est connexe, le groupe $\WGSE$ est un groupe de
Coxeter fini (voir \cite[(5.16)]{ShFC} et \cite[\S~4.2]{ShFCII}). En
g\'en\'eral $\WGSE$ devrait \^etre produit semi-direct d'un groupe de
Coxeter fini par un groupe ab\'elien fini.
Nous notons
$\Phi^+_{\Sigma,\cE}$
l'ensemble des racines positives associ\'ees \`a ce groupe de Coxeter.

Il existe un \'el\'ement unique $w_1$ du groupe de Weyl de $\bG$ tel que
$Z_{\Sigma,\cE}=w_1 W_{\Sigma,\cE}^\bG$ et tel que l'application $\gamma_1\colon 
\WGSE \to \WGSE$ d\'efinie par $\gamma_1(w) =
w_1^{-1}F^{-1}(w) w_1$ envoie tout \'el\'ement de
$\Phi_{\Sigma,\cE}^+$
sur une racine positive. L'isomorphisme $\barQl\WGSE \simeq \End(K)$ donne
lieu \`a un isomorphisme entre les deux diagrammes suivants:
\[
\xymatrix@R=0pt{
\WGSE \ar@/_2pc/[dd]_{F^{-1}} &&
\End(K) \ar@/_2pc/[dd]_{F^*} \\
& Z^\bG_{\Sigma,\cE} \ar[ul]_{w_1^{-1}\cdot} &
& \Hom(K,F^*K) \ar[ul]_{\varphi\circ} \\
W^\bG_{\Sigma,F^*\cE} \ar[ur]_{\cdot w_1}
&& \End(F^*K) \ar[ur]_{\circ\varphi^{-1}}}
\]
Pour $f \in \End(K)$, on a donc $\gamma_1(f) = \varphi \circ F^*(f) \circ
\varphi^{-1}$.

Soit $\tWGSE$ le produit semi-direct de $\WGSE$
par le groupe cyclique engendr\'e par $\gamma_1$. Nous notons
$\Irr(W_{\Sigma,\cE}^\bG)_{\ext}$ l'ensemble des repr\'esentations
irr\'eductibles
de $W_{\Sigma,\cE}^\bG$ qui s'\'etendent en des repr\'esentations de
$\tWGSE$. Pour $E\in\Irr(W_{\Sigma,\cE}^\bG)_{\ext}$, nous choisissons
une extension $\tE$ de $E$ qui est d\'efinie sur $\Q$.

Nous associons \`a tout \'el\'ement $w$ de $\Norm_\bG(\bL)/\bL$ le
sous-groupe
de
L\'evi $\bL^w$ ($F$-stable) de $\bG$ d\'efini comme suit: nous choisissons
un
repr\'esentant $\dot w$ de $w$ dans $\Norm_\bG(\bL)$ ainsi qu'un 
\'el\'ement $g_{\dot w}$
de $\bG$ tel que $g_{\dot w}^{-1}F(g_{\dot w})=F(\dot w)$, et nous posons
\begin{equation} \label{eqn: Lw}
\bL^w:=g_{\dot w}\bL g_{\dot w}^{-1}.
\end{equation}
Soient 
\begin{equation} \label{eqn: SEK}
\Sigma^w=g_{\dot w}\Sigma g_{\dot w}^{-1},\quad \cE^w=\ad(g_{\dot
w}^{-1})^*\cE, \quad K^w=K(\bL^w,\Sigma^w,\cE^w),
\end{equation} et soient
$\varphi_0^w\colon F^*\cE^w\isom\cE^w$ et
$\varphi^w\colon F^* K^w\isom K^w$ les isomorphismes respectivement
induits par $\varphi_0$ et par $\varphi$.

Nous construisons $\varphi_0^{w_1w}\colon F^*\cE^{w_1w}\isom\cE^{w_1w}$ et
$\varphi^{w_1w}\colon F^*K^{w_1w}\isom K^{w_1w}$ au moyen de $\dot
w\circ\ad(\dot w)^*\varphi_0^{w_1}\colon(F\dot w_1\dot w)^*\cE\isom\cE$.
Nous fixons $\varphi_0^{w_1}$, comme il est loisible, de sorte que
$\varphi_0^{w_1}\colon F^* A_0\isom A_0$ co\"\i ncide avec
$\varphi_{A_0}=\varphi_0$.

Le complexe $K^{w_1}$ admet la d\'ecomposition suivante:
\[K^{w_1}=\bigoplus_A A\otimes V_A,\]
o\`u $V_A=\Hom_{\cM G}(A,K^{w_1})$ est une repr\'esentation irr\'eductible
de $W_{\Sigma,\cE}^\bG$.

Pour chaque composante $A$ de $K^{w_1}$, choisissons une structure mixte
$\varphi_A\colon F^*A \isom A$. Ensuite, munissons $V_A$ d'une structure de
$\tWGSE$-module comme suit: pour tout $v \in V_A$, posons
$\gamma_1^{-1}\cdot v = \varphi^{w_1} \circ F^*(v) \varphi_A$. Il est
facile de v\'erifier que cette structure est
bien d\'efinie: pour tout $\theta \in \End(K)$, on a
\begin{multline*}
\gamma_1(\theta) (\gamma_1^{-1}\cdot v)
= (\varphi^{w_1} \circ F^*(\theta) \circ (\varphi^{w_1})^{-1}) \circ 
(\varphi^{w_1} \circ F^*(v) \circ \varphi_A) \\
= \varphi^{w_1} \circ F^*(\theta \circ v) \circ \varphi_A =
\gamma_1^{-1}\cdot(\theta \circ v).
\end{multline*}

Soit $E$ une repr\'esentation irr\'eductible
de $\WGSE$ isomorphe \`a $V_A$. Quitte \`a remplacer $\varphi_A$
par le
produit de celui-ci par une racine de l'unit\'e, nous pouvons supposer
que $V_A$ est isomorphe \`a $\tE$ comme repr\'esentation de
$\tWGSE$, pour
tout $A=A_E$ o\`u $\nu_A$ correspond \`a $\gamma_1^{-1}$ sur $\tilde E$.
Les arguments similaires \`a ceux de \cite[10.4, 10.6]{LCS}, \cite[(2.17),
(5.17)]{ShFC} montrent que
\begin{equation} \label{eqn:KA}
\chi_{K^{w_1w},\varphi^{w_1w}}=\sum_{E\in\Irr(W_{\Sigma,\cE})_{\ext}}
\Tr(\gamma_1 w,\tilde E)\,\chi_{A_E}.
\end{equation}
Il s'ensuit
\begin{equation} \label{eqn:AK}
\chi_{A_E}=|W_{\Sigma,\cE}^\bG|^{-1}\,\sum_{w\in
W_{\Sigma,\cE}^\bG}\Tr(\gamma_1
w,\tilde E)\,\chi_{K^{w_1w},\varphi^{w_1w}}.
\end{equation}

\subsection{La formule du caract\`ere}
\label{subsect:formule}

Pour d\'eterminer la valeur de $\chi_{A_E}$ sur un \'el\'ement $g$ de
$G$, nous sommes donc ramen\'es \`a calculer
$\chi_{K^{w_1w},\varphi^{w_1w}}(g)$. Pour cela, nous \'ecrivons $g=\sigma
v$, o\`u $\sigma$ est semi-simple et $v$ est unipotent et commute \`a
$\sigma$ et nous allons utiliser la formule du caract\`ere qui suit.

Nous notons $\Sigma_\semi$ l'ensemble des parties semi-simples des
\'el\'ements de $\Sigma$, et nous
posons $\bH:=\Cento_\bG(\sigma)$ et $H = \bH^F$. 

Soit $x$ un \'el\'ement de
$G$ tel que $x^{-1}\sigma x\in \Sigma_\semi$. On a alors $\sigma\in x \bL
x^{-1}$.
Posons 
\begin{equation} \label{eqn: Mx}
\bM_x:=(x\bL x^{-1})\cap\bH. 
\end{equation}

Le groupe $\bM_x$ est un sous-groupe de L\'evi
d'un sous-groupe parabolique de $\bH$.

Nous notons $\cO_x$ l'ensemble
des \'el\'ements unipotents $v'$ de $\bH$ tels que $\sigma
v'\in x\Sigma x^{-1}$.
L'ensemble $\cO_x$ est une classe unipotente de $\bM_x$ (\cf
\cite[Proposition~7.11(c)]{LCS}). Soit $\cF_x$
le syst\`eme local sur $\cO_x$, d\'efini comme l'image r\'eciproque de $\cE$
sous
l'application $v\mapsto x^{-1}\sigma vx$ de $\cO_x$ dans $\Sigma$.
Cette application \'etant d\'efinie sur $\F_q$, l'isomorphisme
$\varphi_0\colon F^*\cE\isom\cE$ induit un isomorphisme $\varphi_x\colon
F^*\cF_x\isom\cF_x$.

Soit maintenant $1\boxtimes\cF_x$ l'image r\'eciproque de $\cF_x$ sous
l'application $\Cento(\bM_x)\cO_x\to\cO_x$. La paire
$(\Cento(\bM_x)\cO_x,1\boxtimes\cF_x)$ est une paire cuspidale
$F$-stable sur $\bM_x$. 

\smallskip

Soit $A_x=\IC(\overline{\Cento(\bM_x)\cO_x},1\boxtimes\cF_x)$.
C'est un faisceau caract\`ere cuspidal sur $\bM_x$.
Nous posons $K_x^{\bH}=\ind_{\bM_x}^{\bH}(A_x)$.
La \emph{restriction de $A_x$ \`a la vari\'et\'e unipotente de
$\bM_x$
(et
donc celle de $K_x^{\bH}$ \`a la vari\'et\'e unipotente de
$\bH$)
n'est pas
identiquement nulle}.

La fonction de Green g\'en\'eralis\'ee
$Q^{\bH}_{\bM_x,\cO_x,\cF_x,\varphi_x}$ sur la vari\'et\'e
unipotente de $\bH$ est d\'efinie par (\cf~\cite[(8.3.1)]{LCS}):
\begin{equation} \label{eqn: Green}
Q^{\bH}_{\bM_x,\cO_x,\cF_x,\varphi_x}(v):=
\chi_{K_x^{\bH},\varphi}(v),\qquad
\text{pour tout \'el\'ement unipotent $v$ de $\bH$.}
\end{equation}

On a la formule du caract\`ere suivante (\cite[Theorem 8.5]{LCS}):
\begin{equation} \label{eqn:ffc}
\chi_{K,\varphi_x}(\sigma v)=|H|^{-1}\,|L|^{-1}\,
\sum_{\substack{x\in G\\ x^{-1}\sigma x\in
\Sigma_\semi}}|M_x|\,Q^{\bH}_{\bM_x,\cO_x,\cF_x,
\varphi_x}(v).
\end{equation}

\begin{remark} \label{remark: vers_explicite}
L'\'equation (\ref{eqn:AK}), suivie de (\ref{eqn:ffc}) appliqu\'ee \`a
$\chi_{K^{w_1w},\varphi^{w_1w}}$, puis de (\ref{eqn:KA}) appliqu\'ee \`a
chacun des
$\chi_{K_x^\Gsigma,\varphi_x}$ fournit une certaine expression de
$\chi_{A_E,\varphi_A}$.
Pour obtenir une formule r\'eellement explicite il
faudrait \^etre en mesure de calculer les divers isomorphismes
$\varphi_{?}$ et de d\'ecrire plus pr\'ecis\'ement \[\{x\in G\,:\,
x^{-1}\sigma x\in \Sigma_\semi\}.\] Ce dernier point est trait\'e par Shoji en
\cite[Lemma 4.5]{ShFC}, sous l'hypoth\`ese (v\'erifi\'ee dans les groupes
simples adjoints de type $\BB$, $\CC$ ou $\DD$) que deux \'el\'ements
semi-simples
isol\'es dans un groupe $\bG$ donn\'e\footnote{Un \'el\'ement semi-simple $\sigma$
de $\bG$
est dit \emph{isol\'e} si $\bH$ a m\^eme rang semi-simple que $\bG$.}
et poss\'edant des centralisateurs isomorphes sont conjugu\'es dans $\bG$: si
nous fixons un \'el\'ement $g_0=\sigma_0 v_0$ de $\Sigma$, la
condition $x^{-1}\sigma x\in \Sigma_\semi$ est alors \'equivalente \`a la
condition que
$\bM_x$ soit un sous-groupe de L\'evi de $\bH$ conjugu\'e sous $\bG$
\`a $\Cento_\bG(\sigma_0)$.
\end{remark}

Soit $w\in \WGL$. Nous noterons $k_{\bL,\Sigma,\cE}(w)=k(w)$
la fonction trace du complexe tordu $K^w$. 
L'\'equation~(\ref{eqn:ffc}) appliqu\'ee \`a $K^w$ donne:
\begin{equation} \label{eqn:k(w)}
k(w)(\sigma v)=|H|^{-1}\,|L^{w}|^{-1}\,
\sum_{\substack{x\in G\\ x^{-1}\sigma x\in
\Sigma_\semi^w}}|M^w_x|\,Q^{\bH}_{\bM^w_x,\cO_x^w,\cF_x^w,\varphi_x^w}(v),
\end{equation}
o\`u
\[\bM_x^w=\Cent_\bH(x\bT_1^wx^{-1})=x\Cento_{\bL^w}(x^{-1}\sigma x) x^{-1}.
\]
\section{Notations}
\label{sect:notations}

\subsection{Sous-groupes de L\'evi dans $\bH$}

On garde toutes les notations du paragrphe pr\'ec\'edent: $\bG$ est un groupe r\'eductif; $\bL$ est un sous-groupe de L\'evi poss\'edant une paire cuspidale $(\Sigma,\cE)$; $\Sigma_\semi$ est l'ensemble des parties semi-simples des \'el\'ements de $\Sigma$; et $g = \sigma v$ est l'\'el\'ement en lequel on veut calculer la valeur de $\chi_{A_E}$. Consid\'erons l'ensemble
\begin{equation} \label{eqn:cM}
\cM = \{ \bM_x : x^{-1} \sigma x \in \Sigma_\semi \},
\quad\text{o\`u}\qquad
\bM_x = \Cent_\bH(x\bT_1 x^{-1}) = x\Cento_\bL(x^{-1}\sigma x) x^{-1}.
\end{equation}
(On n'exige pas que $x \in G$ ici). Ainsi que nous l'avons d\'ej\`a remarqu\'e, 
tout $\bM_x$ est un sous-groupe de L\'evi du groupe r\'eductif $\bH$. Rappelons 
que les sous-groupes de L\'evi d'un groupe r\'eductif \`a conjugaison pr\`es sont 
param\'etr\'es par les orbites du groupe de Weyl sur l'ensemble des 
sous-graphes du graphe de Dynkin. En particulier, il y a, \`a conjugaison 
pr\`es, un nombre fini de sous-groupes de L\'evi. Par cons\'equent, $\cM$ est
partitionn\'e  en un nombre fini de classes d'\'equivalence
\begin{equation} \label{eqn:cMj}
\cM = \cM_1 \sqcup \cdots \sqcup \cM_r,
\end{equation}
o\`u $\bM_x$ et $\bM_y$ sont dits \'equivalents s'ils sont conjugu\'es sous 
$\bH$. (Nous verrons \`a la Section~\ref{sect:sousgp} que si $\bG$ est 
quasi-simple, alors $r \le 2$, et en fait $r = 1$ dans la plupart des cas).

\begin{remark}\label{rmk:m-unique}
Il est \`a noter que tous les membres de $\cM$ sont conjugu\'es sous $\bG$, et donc isomorphes. En effet, il est clair que chaque $\bM_x$ est conjugu\'e sous $\bG$ au centralisateur dans $\bL$ d'un \'el\'ement de $\Sigma_\semi$. Mais il s'ensuit de la d\'efinition de $\Sigma$ que tous les \'el\'ements de $\Sigma_\semi$ sont conjugu\'es \`a multiplication par un \'el\'ement central pr\`es, et donc leurs centralisateurs sont conjugu\'es dans $\bL$.
\end{remark}

Pour tout $j \in \{1, \ldots, r\}$, nous fixons un \'el\'ement $g_{a_j} \in \bG$ 
tel que $\bM_{g_{a_j}} \in \cM_j$. Posons
\begin{equation} \label{eqn:Mj}
\bT^\bH_j := g_{a_j} \bT_1 g_{a_j}^{-1}
\qquad\text{et}\qquad
\bM_j := \bM_{g_{a_j}} = \Cent_\bH(\bT^\bH_j).
\end{equation}
Ensuite, soit $a_j$ l'\'el\'ement de $\WGL$ tel que $F(a_j)$ soit l'image de $g_{a_j}^{-1} F(g_{a_j})$ dans $\WGL$.
Le groupe $\bM_j$ est un sous-groupe de
L\'evi de $\bH$. Nous posons
\begin{equation} \label{eqn:WHM}
\WHMj :=\Norm_{\bH}(\bM_j)/\bM_j.
\end{equation}

\begin{exm} \label{ex:M}
Supposons qu'il existe un \'el\'ement $x_1$ de $G$ tel que 
$x_1^{-1}\sigma x_1\in\Sigma_\semi^w$. Pour l'un des $g_{a_j}$, on peut prendre $g_{a_j} = x_1 g_{\dot w}$. On aura alors
$\bT^\bH_j =x_1 \bT^w_1 x_1^{-1}$ et 
$\bM_j =x_1\Cento_{\bL^w}(x_1^{-1}\sigma x_1)x_1^{-1}\subset \bH$.
Dans ce cas $a_j=w$.
\end{exm}
  
On remarque aussi que $\Norm_\bH(\bT^\bH_j)=\Norm_\bH(\bM_j)$ et
$\Norm_\bG(\bT_1)=\Norm_\bG(\bL)$. 

De l'application 
\begin{equation} \label{eqn:plongtI}
\tilde\iota_j\colon\Norm_\bH(\bM_j)\to\Norm_\bG(\bL),\quad h\mapsto
g_{a_j}^{-1}hg_{a_j}\end{equation}
se d\'eduit un plongement 
\begin{equation} \label{eqn: plongtII}
\iota_j\colon \WHMj\hookrightarrow \WGL.
\end{equation}

Bien que le morphisme de Frobenius agisse et sur $\WGL$ et sur tous les $\WHMj$, 
il est \`a noter que le plongement $\iota_j$ n'est en g\'en\'eral 
pas \'equivariant pour ces deux actions. Nous gardons la notation $F\colon \WGL \to \WGL$ pour l'automorphisme de $\WGL$ induit par le Frobenius, et note $\eta_j\colon\WHMj \to \WHMj$ l'automorphisme \emph{inverse} \`a celui qui est induit par le Frobenius sur $\WHMj$. Il s'ensuit des d\'efinitions de $\tilde\iota_j$ et de $a_j$ que
\begin{multline*}
\tilde\iota_j(F(h)) = g_{a_j}^{-1}F(h) g_{a_j} \\
= g_{a_j}^{-1}F(g_{a_j}) F(g_{a_j}^{-1}h g_{a_j}) F(g_{a_j})^{-1}g_{a_j}
= F(a_j) F(\tilde\iota_j(h)) F(a_j)^{-1}
\end{multline*}
Rempla\c cons $F(h)$ par son image $w$ dans $\WHMj$ et $\tilde\iota_j$ par $\iota_j$. On obtient l'\'egalit\'e $\iota_j(w) = F(a_j) F(\iota_j(\eta_j(w))) F(a_j)^{-1}$, ou, autrement dit, 
\begin{equation}\label{eqn:eta-F}
a_j \iota_j(\eta_j(w)) = F^{-1}(\iota_j(w)) a_j.
\end{equation}
D\'esormais, nous identifions $\WHMj$ avec son image par $\iota_j$. Si $w \in \WHMj$, il n'y a donc aucune ambigu\"\i t\'e dans la notation $F(w)$: c'est l'image sous le morphisme de Frobenius de $w$ en tant qu'\'el\'ement de $\WGL$.

\subsection{Actions de Frobenius et Doubles Classes dans $\WGL$}

Il s'ensuit de~\eqref{eqn:eta-F} que $\WHMj$ agit sur $a_j \WHMj$ par $F^{-1}$-conjugaison; \ie pour tout $w \in \WHMj$, on a $F^{-1}(w)(a_j \WHMj)w^{-1} = a_j \WHM_j$. De plus, les orbites de cette action sont en bijection avec les classes de $\eta_j$-conjugaison dans $\WHMj$, via la bijection \'evidente $\WHMj \leftrightarrow a_j \WHMj$ donn\'ee par $w \leftrightsquigarrow a_jw$. (Rappelons que deux \'el\'ements $u,v \in \WHMj$ sont dits \emph{$\eta_j$-conjugu\'es} s'il existe un $w \in \WHMj$ tel que $\eta_j(w)u w^{-1} = v$).

Cette situation est parall\`ele \`a celle de $\WGSE$ et $w_1\WGSE$. \`A savoir, $\WGSE$ agit sur $w_1\WGSE$ par $F^{-1}$-conjugaison, et les orbites de cette action sont en bijection avec les classes de $\gamma_1$-conjugaison dans $\WGSE$.

Nous notons $\sim_{\eta_j}$, $\sim_{\gamma_1}$, et $\sim_{F^{-1}}$ les relations de $\eta_j$-, $\gamma_1$-, et $F^{-1}$-conjugaison, respectivement. Si $w \in \WGSE$, son stabilisateur sous la $\gamma_1$-conjugaison sera not\'e $\Cent_{\gamma_1}(w)$. Si $w \in \WHMj$ (resp.~$w \in \WGL$), alors $\Cent_{\eta_j}(w)$ (resp.~$\Cent_{F^{-1}}(w)$) est d\'efini de mani\`ere semblable.

Ensuite, nous notons $\tWHMj$ le produit semi-direct de $\WHMj$ par le groupe cyclique engendr\'e par $\eta_j$. $\Irr(\WHMj)_\ext$ d\'esigne l'ensemble des repr\'esentations irr\'eductibles qui admettent une action de $\tWHMj$. Pour tout $E' \in \Irr(\WHMj)_\ext$, nous fixons une fois pour toutes une extension $\tE'$ en repr\'esentation irr\'educible de $\tWHMj$.

Les calculs de la Section~\ref{sect:princ} entra\^\i neront une comparaison de l'action de $\WHMj$ sur $a_j \WHMj$ et celle de $\WGSE$ sur $w_1 \WGSE$. 
Nous posons
\[
\cN_j := \WGSE \backslash \WGL / \WHMj.
\]

Dans chaque double classe $\unu \in \cN_j$, choisissons, une fois pour toutes, un repr\'esentant $w_\unu$.  Nous posons 
\begin{equation} \label{Wunu}
W(\unu) = w_1 \WGSE \cap F^{-1}(w_\unu) (a_j \WHM_j) w_\unu^{-1},
\end{equation}
et nous d\'efinissons deux plongements comme suit:
\[
\begin{aligned}
\lambda\colon& W(\unu) \to \WHMj, \qquad& \lambda(\uw) &= \eta_j^{-1}(a_j^{-1}
F^{-1}(w_\unu^{-1}) \uw w_\unu) \\
\kappa\colon& W(\unu) \to \WGSE, \qquad& \kappa(\uw) &= \gamma_1^{-1}(w_1^{-1}
\uw)
\end{aligned}
\]
Nous remarquons que si l'on remplace $w_\unu$ par $w_\unu v$, o\`u $v \in \WHMj$, alors $W(\unu)$ ne change pas, car $a_j \WHM_j$ est stable sous $F^{-1}$-conjugaison par $v$. D'autre part, si l'on remplace $w_\unu$ par $y w_\unu$, o\`u $y \in \WGSE$, alors $W(\unu)$ est remplac\'e par son $F^{-1}$-conjugu\'e 
$F^{-1}(y)W(\unu)y^{-1} \subset w_1 \WGSE$.

\begin{lem}\label{lem:wunu-card}
Soit $t \in \unu$, et posons $P(t) = \{(y,z) \in \WGSE \times \WHMj : y
w_\unu z = w_1^{-1} F^{-1}(t) a_j \}$.  Alors $|P(t)| = |W(\unu)|$.
\end{lem}
\begin{proof}
Montrons le lemme d'abord dans le cas o\`u $t = w_\unu$. Si $y w_\unu z =
w_1^{-1} F^{-1}(w_\unu) a_j$, alors $w_1 y = F^{-1}(w_\unu) a_jz^{-1}
w_\unu^{-1}$, et donc $w_1 y \in W(\unu)$. L'application $\phi\colon P(w_\unu) \to
W(\unu)$ d\'efinie par $(y,z) \mapsto w_1 y$ est injective puisque $z$ est
d\'etermin\'e par $y$. D'autre part, cette application est aussi surjective: si $w \in W(\unu)$, alors posons $y = w_1^{-1}w \in \WGSE$ et $z = (a_j^{-1}F^{-1}(w_\unu^{-1}) w w_\unu)^{-1} \in \WHMj$. Il est clair que $y w_\unu z = w_1^{-1} f^{-1}(w_\unu) a_j$, et donc $(y,z) \in P(w_\unu)$ et $\phi(y,z) = w$.

Ensuite, si $t \ne w_\unu$, posons $W'(\unu) = w_1 \WGSE \cap F^{-1}(t) (a_j 
\WHMj) t^{-1}$. D'une part, l'argument du paragraphe pr\'ec\'edent montre que 
$|P(t)| = |W'(\unu)|$, et d'autre part, les remarques qui pr\'ec\`edent le lemme 
montrent que $|W(\unu)| = |W'(\unu)|$.
\end{proof}

Nous introduisons maintenant une famille d'accouplements, param\'etr\'ee par les $\cN_j$, qui relient les ensembles
$\Irr(\WGSE)_\ext$ et $\Irr(\WHMj)_\ext$.  Soit
$\unu \in
\cN_j$.  Si $E \in \Irr(\WGSE)_\ext$ et $E' \in
\Irr(\WHMj)_\ext$, on pose
\[
\langle E, E' \rangle_\unu 
:= \frac{1}{|W(\unu)|} \sum_{\uw \in W(\unu)} \Tr(\gamma_1 \kappa(\uw), \tilde E)
\Tr(\eta_j \lambda(\uw), \tilde E').
\]

\subsection{Fonctions de Green}

Pour tout $j$ et tout $w \in \WHMj$, nous choisissons un repr\'esentant $\dot w \in \Norm_\bH(\bM_j)$ ainsi qu'un \'el\'ement $h_{\dot w} \in \bH$ tel que $h_{\dot w}^{-1} F(h_{\dot w}) = F(\dot w)$. Ensuite, posons
\[
\bM_j^w = h_{\dot w} \bM_j h_{\dot w}^{-1}.
\]
Les $\bM_j^w$ sont tous $F$-stables. $\bM_j^w$ et $\bM_j^{w'}$ sont conjugu\'es sous $H$ si et seulement si $w$ et $w'$ sont $\eta_j$-conjugu\'es dans $\WHMj$.

R\'evisons la construction de la fonction de Green g\'en\'eralis\'ee 
$Q^\bH_{\bM_x,\cO_x,\cF_x,\varphi_x}$. Nous allons d'abord v\'erifier que sa 
d\'efinition ne d\'epend pas de $x$, mais seulement du groupe $\bM_x$. 
(C'est-\`a-dire, si $\bM_x = \bM_y$, alors $Q^\bH_{\bM_x, \cO_x, \cF_x, 
\varphi_x} = Q^\bH_{\bM_y, \cO_y, \cF_y, \varphi_y}$). 
De plus, si $\bM_x$ et $\bM_y$ sont conjugu\'es sous $H$, alors leurs fonctions 
de Green respectives sont \'egales.

\begin{remark} \label{rem:classesconj}
Nous supposons que les groupes $\bM_x$ et $\bM_y$ \'egaux et que, 
soit $\bH$, soit $\bL$, est \'egal \`a $\bG$.
L'\'egalit\'e des groupes 
$\bM_x$ et $\bM_y$ est alors \'equivalente \`a 
celle de $x\bL x^{-1}$ et de $y\bL y^{-1}$ (en effet, si $\bH=\bG$, nous
avons
$\bM_x=x\bL x^{-1}$, et si $\bL=\bG$, nous avons $x\bL x^{-1}=\bG$, et, par
cons\'equent,
$\bM_x=\bH$). Il existe alors 
$n\in\Norm_{\bG}(\bL)$ tel que $y=nx$ et l'orbite unipotente $\cO_y$ s'\'ecrit 
\[\cO_y=\left\{v''\in\bG\,:\,\sigma v''\in nx\Sigma x^{-1}n^{-1}\right\}
=\left\{v''\in\bG\,:\,\sigma(n^{-1}v''n)\in x\Sigma x^{-1}\right\},\]
puisque $\sigma$ est central dans $\bG$, autrement dit, 
$\cO_y=n\cO_x n^{-1}$. 
\end{remark}

\begin{lem} \label{lem:indepGreenfunct}
Si $\bM_x = \bM_y$, alors $Q^\bH_{\bM_x, \cO_x, \cF_x, 
\varphi_x} = Q^\bH_{\bM_y, \cO_y, \cF_y, \varphi_y}$.
\end{lem}\begin{proof}

Nous supposons les groupes $\bM_x$ et $\bM_y$ \'egaux. 

$\bullet$ {\em \'Egalit\'e de $\cO_x$ et de $\cO_y$.}
Si $\bM_x=\bM_y=\bG$, on a en particulier $\bH=\bG$, et donc il r\'esulte de la remarque~\ref{rem:classesconj} que $\cO_x=\cO_y$.

Nous supposons dor\'enavant $\bM_x=\bM_y\neq\bG$. Lorsque 
$\bG$ est un groupe de type exceptionnel,  la
Table~\ref{tbl:levi} montre que les seuls cas \`a 
consid\'erer sont les suivants:
\begin{enumerate}
\item
$\bG=\EE_7^\ad$ ou $\bG=\EE_8$ et $\bM_x=x\bL x^{-1}=\EE_6$ ($p=3$): d'apr\`es 
\cite[Proposition~20.3.(a)]{LCS} la classe unipotente r\'eguli\`ere de $\bM_x$ 
est la seule \`a porter un syst\`eme local cuspidal,
\item
$\bG=\EE_8$, $\bM_x=x\bL x^{-1}=\EE_7$ ($p=2$): d'apr\`es 
\cite[Proposition~20.3.(c)]{LCS} la classe unipotente r\'eguli\`ere de
$\bM_x$
est la seule \`a porter un syst\`eme local cuspidal,
\item
$\bG=\bL=\EE_8$, $\bM_x=\EE_7\times\AA_1$ ($p\ne 2$): il n'y a pas de
syst\`eme local cuspidal \`a support unipotent dans $\EE_7$ lorsque $p\ne 2$,
donc ce cas ne se produit pas,
\item
$\bG=\bL=\EE_8$, $\bM_x=\EE_6\times\AA_2$ ($p\ne 2$):
\begin{itemize}
\item[(a)]
si $p\ne 3$, d'apr\`es \cite[Proposition~20.3]{LCS} la classe unipotente
$\EE_6(a_3)$ est la seule \`a porter des syst\`eme locaux cuspidaux,
\item[(b)]
si $p=3$, la classe unipotente
r\'eguli\`ere de $\EE_6$ est la seule \`a porter des syst\`eme locaux
cuspidaux.
\end{itemize}
\end{enumerate}
Dans chacun des cas \'enum\'er\'es ci-dessus, il y a donc au plus une
classe unipotente portant des syst\`eme locaux cuspidaux, donc
n\'ecessairement $\cO_x=\cO_y$.

$\bullet$ {\em Isomorphie de $\cF_x$ et de $\cF_y$.}
Dans chacun des cas (1), (2) et (4), 
la remarque~\ref{rem:classesconj} s'applique et montre qu'il existe
$n\in\Norm_{\bG}(\bL)$ tel que $y=nx$. 
D'apr\`es \cite[(8.3.2)]{LCS}, pour d\'efinir la fonction de Green 
g\'en\'eralis\'ee $Q^\bH_{\bM_x,\cO_x,\cF_x,\varphi_x}$, nous aurions pu
remplacer $\cE$ par un syst\`eme local $\cE_1$ d\'efini comme 
l'image r\'eciproque sous l'application naturelle 
$\bL\to \bL/\bL_{\der}\times \bL/\bT_1$, de
$\barQl\boxtimes\cE'$, o\`u $\cE'$ est un syst\`eme local irr\'eductible
$\bL$-\'equivariant sur $\cO$.
Nous savons, d'apr\`es \cite[Theorem~9.2]{LIC}, que le 
groupe $W^\bG_{\Sigma,\cE_1}$ est isomorphe \`a
$\Norm_{\bG}(\bL)/\bL$.
L'\'el\'ement $n$ v\'erifie donc
$(\ad(n))^*\cE_1\simeq \cE_1$. Nous avons vu que $\cO_x=\cO_y$. Il
s'ensuit que $\cF_x\simeq\cF_y$. 

$\bullet$ {\em Ind\'ependance de la structure mixte.}
L'argument qui suit est inspir\'e de celui utilis\'e par Shoji dans la preuve 
de \cite[Lemme~4.5]{Sh}. Nous notons encore $\varphi_0$ l'isomorphisme
$\varphi_0\colon F^*\cE_1\isom\cE_1$ induit par $\varphi'\colon 
F^*\cE'\isom \cE'$.
Par d\'efinition de l'isomorphisme $\varphi_x\colon F^*\cF_x\isom\cF_x$, on a 
\[\Tr(\varphi_x,(\cF_x)_{v'})=\Tr(\varphi_0, (\cE_1)_{x^{-1}\sigma v'x}),\quad
\text{pour tout $v'\in\cO_x^F$.}\]
Notons $\bar n$ l'image de $n$ par la projection $\bL\to\bL/\bT_1$.
L'action par conjugaison de $\bar n$ sur $\bL/\bT_1$ induit un
isomorphisme $(\ad \bar n)^*\cE'\isom \cE'$ compatible avec l'isomorphisme
$\varphi'\colon F^*\cE'\isom\cE'$. Il en r\'esulte que
$\Tr(\varphi',\cE'_{p(x^{-1}\sigma v'x)})=
\Tr(\varphi',\cE'_{p(y^{-1}\sigma v'y)})$, o\`u $p$ est la projection
naturelle de $\Sigma$ sur $\Sigma/\bT_1$. L'\'egalit\'e de
$\Tr(\varphi_x,(\cF_x)_{v'})$ et de $\Tr(\varphi_y,(\cF_y)_{v'})$
s'ensuit. La structure mixte ne d\'epend donc que de $\bM_x$.
\end{proof} 

Pour tout $w \in \WHMj$, posons
\begin{equation} \label{eqn:Qw}
Q_w := \text{la fonction de Green associ\'ee \`a $\bM_j^w$.}
\end{equation}
D'apr\`es le paragraphe pr\'ec\'edent, $Q_w$ est bien d\'efini, et
\[
Q_w = Q_{w'}
\qquad\text{si $w$ et $w'$ sont $\eta_j$-conjugu\'es.}
\]

Enfin, pour tout $E' \in \Irr(\WHMj)_\ext$, on pose
\[
Q_{E'}(v) := \frac{1}{|\WHMj|} \sum_{w' \in \WHMj} \Tr(\eta_j w', \tilde E') Q_{w'}(v).
\]
Il s'ensuit que
\[
Q_{w'}(v) = \sum_{E' \in \Irr(\WHMj)_\ext} \Tr(\eta_j w', \tilde E') Q_{E'}(v).
\]

\section{La Formule du Caract\`ere}
\label{sect:princ}

Dans cette section, nous \'etablissons le r\'esultat principal en trois \'etapes. La premi\`ere \'etape (le Lemme~\ref{lem:k-mw}) consiste \`a modifier la formule~\eqref{eqn:k(w)} de mani\`ere qu'il n'y reste aucune mention d'\'el\'ements de $G$. Dans la seconde \'etape (le Lemme~\ref{lem:k-w}), on fait intervenir les ensembles de doubles classes $\cN_j$ dans la formule. Cela permet enfin d'\'ecrire au Th\'eor\`eme~\ref{thm:princ} la formule cherch\'ee pour $\chi_{A_E}$.

\begin{lem}[{\cf \cite[Proposition~2.16]{MW}}]\label{lem:k-mw}
Pour tout $w \in \WGSE$, on a
\[
k(w_1w)(\sigma v) = \sum_{j=1}^r \frac{|\Cent_{F^{-1}}(w_1w)|}{|\WHMj|} 
\sum_{\substack{w' \in
\WHMj \\ a_jw' \sim_{F^{-1}} w_1w}} Q_{w'}
\]
\end{lem}
\begin{proof}
D'apr\`es~\eqref{eqn:k(w)}, on a
\begin{equation} \label{eqn:k(w1w)}
k(w_1w)(\sigma v) = |H|^{-1} |L^{w_1w}|^{-1} \sum_{\substack{x \in G\\
x^{-1}\sigma x \in \Sigma^{w_1w}_\semi}} |M^{w_1w}_x| Q^\bH_{\bM^{w_1w}_x, 
\cO^{w_1w}_x, \cF^{w_1w}_x, \varphi^{w_1w}_x} (v).
\end{equation}
Soit $x \in G$ tel que $x^{-1}\sigma x \in \Sigma^{w_1w}_\semi$. Il s'ensuit que
$g_{\dot{(w_1w)}}^{-1}x^{-1}\sigma xg_{\dot{(w_1w)}}\in\Sigma_\semi$.
\'Evidemment, on a
\[
\bM^{w_1w}_x = \Cent_\bH(x\bT_{w_1w}x^{-1}) = \Cent_\bH(x g_{\dot{(w_1w)}} \bT_1 g_{\dot{(w_1w)}}^{-1} x^{-1}).
\]
En particulier, $\bM^{w_1w}_x$ appartient \`a l'ensemble $\cM$, et donc \`a un 
certain $\cM_{j}$, o\`u $j=j(w,x) \in \{1, \ldots, r\}$. Ensuite, 
$\bM^{w_1w}_x$, \'etant $F$-stable, est conjugu\'e sous $H$ \`a un $\bM_{j}^{w'}$, o\`u $w' \in \WHMj$ est unique \`a $\eta_j$-conjugaison pr\`es. Cette classe de $\eta_j$-conjugaison sera not\'ee $C(w,x)$. Pour tout $w' \in C(w,x)$, on a
\[
Q_{w'} = Q^\bH_{\bM^{w_1w}_x, \cO^{w_1w}_x, \cF^{w_1w}_x, \varphi^{w_1w}_x}.
\]
De plus, si $w' \in C(w,x)$, le fait que
\[
\Cento(\bM_j^{w'}) = g_{\dot w'} g_{a_j} \bT_1 g_{a_j}^{-1} g_{\dot w'}^{-1}
\qquad\text{et}\qquad
\Cento(\bM^{w_1w}_x) = x g_{\dot{(w_1w)}} \bT_1 g_{\dot{(w_1w)}}^{-1} x^{-1}
\]
sont conjugu\'es sous $G$ implique que les images dans $\WGL$ des deux 
\'el\'ements suivants sont $F^{-1}$-conjugu\'ees:
\begin{align*}
g_{a_j}^{-1} g_{\dot w'}^{-1} F(g_{\dot w'} g_{a_j}) &= g_{a_j}^{-1} F(\dot w') g_{a_j} g_{a_j}^{-1} F(g_{a_j}) = \tilde\iota(F(\dot w')) g_{a_j}^{-1} F(g_{a_j}), \\
g_{\dot{(w_1w)}}^{-1} x^{-1} F(x g_{\dot{(w_1w)}}) &= F(\dot{(w_1w)}).
\end{align*}
Leurs images dans $\WGL$ sont $\eta_j^{-1}(w')F(a_j)$ et $F(w_1w)$, respectivement. De la d\'efinition de $\eta_j$ se d\'eduit la formule
$\eta_j^{-1}(w') = F(a_j w' a_j^{-1})$. On sait donc que $F(a_j w') \sim_{F^{-1}} F(w_1w)$. Il est clair que $F$ respecte les classes $F^{-1}$-conjugaison, et la condition pr\'ec\'edente \'equivaut \`a ce que
\[
a_j w' \sim_{F^{-1}} w_1w.
\]

Pour chaque classe de $\eta_j$-conjugaison $C \subset \WHMj$, posons
\[
Y(C) = \{ x \in G : \text{$x^{-1}\sigma x \in \Sigma_\semi^{w_1w}$, $j(w,x) = j$ 
et $C(w,x) = C$} \}.
\]
L'ensemble $Y(C)$ est vide sauf si un (et donc tout) membre de $a_jC$ est 
$F^{-1}$-conjugu\'e \`a $w_1w$. Choisissons un repr\'esentant $w'_C$ de
chaque classe. Le cardinal $|M_j^{w'_C}|$ et la fonction de Green $Q_{w'_C}$ sont tous deux ind\'ependants du choix de $w'_C$. La formule~(\ref{eqn:k(w1w)})
devient donc:
\[
k(w_1w)(\sigma v) = |H|^{-1} |L^{w_1w}|^{-1} \sum_{j=1}^r 
\sum_{\substack{C \subset
\WHMj \\ a_jw'_C \sim_{F^{-1}} w_1w}} \sum_{x \in Y(C)} |M^{a_jw'_C}| Q_{w'_C}(v).
\]
L'\'el\'ement $x$ ne joue plus aucun r\^ole dans la formule; on peut
remplacer la deuxi\`eme somme ci-dessus par le cardinal de $Y(C)$. Ensuite,
on peut remplacer la somme sur certaines classes de
$\eta_j$-conjugaison par une somme sur certains \'el\'ements $w' \in
\WHMj$, pourvu qu'on divise aussi chaque terme par le cardinal de la classe
correspondante. Ce cardinal est \'egal \`a $|\WHMj|/|\Cent_{\eta_j}(w')|$, et
l'on obtient donc
\[
k(w_1w)(\sigma v) = |H|^{-1} |L^{w_1w}|^{-1} \sum_{j=1}^r 
\sum_{\substack{w' \in \WHMj\\ a_jw' \sim_{F^{-1}} w_1 w}}
\frac{|Y(w')||M^{a_jw'_C}||\Cent_{\eta_j}(w')|}{|\WHMj|} Q_{w'}(v).
\]
Ici, $Y(w')$ d\'esigne l'ensemble $Y(C)$, o\`u $C$ est la classe de
$\eta_j$-conjugaison de $w'$. Enfin, l'argument de~\cite[p.~510]{MW}
montre que
\[
|Y(w')| = \frac{|H| |L^{w_1w}|
|\Cent_{F^{-1}}(w_1w)|}{|M^{a_jw'_C}||\Cent_{\eta_j}(w')|}.
\]
La formule cherch\'ee s'ensuit. 
\end{proof}

\begin{lem}[{\cf \cite[Lemme~7.1]{W}}]\label{lem:k-w}
Pour tout $w \in \WGSE$, on a
\[
k(w_1w)(\sigma v) = \sum_{j=1}^r
|\Cent_{\eta_j}(w)| \sum_{\unu \in \cN_j} \frac{1}{|W(\unu)|}
\sum_{\substack{\uw \in W(\unu) \\ \kappa(\uw) \sim_{\gamma_1} w}} Q_{\lambda(\uw)}.
\]
\end{lem}
\begin{proof}
Il est clair que pour tout $w'$ figurant dans la somme du
Lemme~\ref{lem:k-mw}, le cardinal de l'ensemble des $t \in \WGL$ tels que
$F^{-1}(t)a_jw' t^{-1} = w_1 w$ \'egale celui de $\Cent_{F^-1}(w_1w)$. Donc
\[
k(w_1w)(\sigma v) = \sum_{j=1}^r
\frac{1}{|\WHMj|} \sum_{\substack{w' \in \WHMj\\ a_jw' \sim_{F^{-1}} w_1
w}} \ \sum_{\substack{t \in \WGL\\ F(t) a_jw' t^{-1} = w_1w}} Q_{w'}
\]
Pour tout $t$, soit $\unu(t)$ la double classe $\WGSE t \WHMj$. D'apr\`es le
Lemme~\ref{lem:wunu-card}, on a
\begin{multline*}
k(w_1w)(\sigma v)\\
\begin{aligned}
&= \sum_{j=1}^r
\frac{1}{|\WHMj|} \sum_{\substack{w' \in \WHMj\\ a_jw' \sim_{F^{-1}} w_1
w}} \ \sum_{\substack{t \in \WGL\\ F^{-1}(t) a_jw' t^{-1} = w_1w}} 
\frac{1}{|W(\unu(t))|} \sum_{\substack{y \in \WGSE,\ z \in \WHMj\\ y
w_\unu z = w_1^{-1} F^{-1}(t) a_j}} Q_{w'} \\
&= \sum_{j=1}^r
\frac{1}{|\WHMj|} \sum_{\unu \in \cN_j} \frac{1}{|W(\unu)|}
\sum_{\substack{w' \in \WHMj,\ t \in \WGL\\ \unu(t) = \unu\\ F^{-1}(t) a_jw'
t^{-1} = w_1 w}} \ \sum_{\substack{y \in \WGSE,\ z \in \WHMj\\ y
w_\unu z = w_1^{-1} F^{-1}(t) a_j}} Q_{w'}.
\end{aligned}
\end{multline*}
La condition $y w_\unu z = w^{-1} F^{-1}(t) a_j$ implique que $t^{-1} = 
F(a_j^{-1} w_\unu^{-1} y^{-1} w_1^{-1})$, et la condition $F^{-1}(t) a_j w' t^{-1} = w_1 w$ \'equivaut \`a ce que $y w_\unu z w' F(a_j z^{-1} w_\unu^{-1} y^{-1} w_1^{-1}) = w$. La formule ci-dessus devient donc:
\begin{equation}\label{eqn:lem-kw1}
k(w_1w)(\sigma v) = \sum_{j=1}^r
\frac{1}{|\WHMj|} \sum_{\unu \in \cN_j} \frac{1}{|W(\unu)|}
\sum_{\substack{w' \in \WHMj,\ y \in \WGSE,\ z \in \WHMj\\
y w_\unu z w' F(a_j z^{-1} w_\unu^{-1} y^{-1} w_1^{-1}) = w}} Q_{w'}.
\end{equation}
La condition dans la troisi\`eme somme peut \'egalement s'\'ecrire sous 
de nombreuses formes \'equivalentes, dont:
\begin{align*}
y w_\unu z w' F(a_j z^{-1} w_\unu^{-1} y^{-1} w_1^{-1}) &= w \\
w_\unu zw' F(a_j z^{-1} w_\unu^{-1}) &= y^{-1} w F(w_1 y) \\
F^{-1}(w_\unu) F^{-1}(zw') a_j z^{-1} w_{\unu}^{-1} &= F^{-1}(y^{-1}w) w_1 y \\
F^{-1}(w_\unu) a_j\eta_j(zw')z^{-1} w_{\unu}^{-1} &= w_1 \gamma_1(y^{-1} w) y \\
F^{-1}(w_\unu) a_j\eta_j(zw' \eta_j^{-1}(z^{-1}) ) w_\unu^{-1} &= w_1 \gamma_1( y^{-1} w \gamma_1^{-1}(y))
\end{align*}
La derni\`ere de ces possibilit\'es \'equivaut \`a l'\'enonc\'e suivant:
\[
\text{Il existe un $\uw \in W(\unu)$ tel que
$\kappa(\uw) = y^{-1} w \gamma_1^{-1}(y)$
et
$\lambda(\uw) = z w' \eta_j^{-1}(z^{-1})$.}
\]
On peut maintenant \'ecrire la formule comme suit:
\[
k(w_1w)(\sigma v) = \sum_{j=1}^r
\frac{1}{|\WHMj|} \sum_{\unu \in \cN_j} \frac{1}{|W(\unu)|}
\sum_{\substack{w' \in \WHMj,\ \uw \in W(\unu),\ y \in \WGSE,\ z \in \WHMj\\
\kappa(\uw) = y^{-1} w \gamma_1^{-1}(y),\ \lambda(\uw) = z w' \eta_j^{-1}(z^{-1})}}
 Q_{w'}.
\]
Il est clair qu'on peut remplacer l'\'egalit\'e $\kappa(\uw) = y^{-1} w \gamma_1^{-1}(y)$ par la condition $\kappa(\uw) \sim_{\gamma_1} w$, pourvu qu'on introduise en m\^eme temps la multiplicit\'e $|\Cent_{\gamma_1}(w)|$, et de m\^eme pour la condition sur $w'$:
\[
k(w_1w)(\sigma v) = \sum_{j=1}^r
\frac{1}{|\WHMj|} \sum_{\unu \in \cN_j} \frac{1}{|W(\unu)|}
\sum_{\substack{w' \in \WHMj,\ \uw \in W(\unu)\\
\kappa(\uw) \sim_{\gamma_1} w,\ \lambda(\uw) \sim_{\eta_j} w'}}
|\Cent_{\gamma_1}(w)| |\Cent_{\eta_j}(w')| Q_{w'}.
\]
Ensuite, puisque $Q_{w'}$ ne d\'epend que la classe de
$\eta_j$-conjugaison de $w'$, on peut le remplacer par
$Q_{\lambda(\uw)}$. On peut \'egalement remplacer
$|\Cent_{\eta_j}(w')|$ par $|\Cent_{\eta_j}(\lambda(\uw))|$, et on
obtient ainsi une formule dans laquelle $w'$ ne joue plus de grand r\^ole:
\[
k(w_1w)(\sigma v) = \sum_{j=1}^r
\frac{|\Cent_{\gamma_1}(w)|}{|\WHMj|}
\sum_{\unu \in \cN_j} \frac{1}{|W(\unu)|}
\sum_{\substack{\uw \in W(\unu)\\
\kappa(\uw) \sim_{\gamma_1} w}}\ 
\sum_{\substack{w' \in \WHMj\\ \lambda(\uw) \sim_{\eta_j} w'}}
 |\Cent_{\eta_j}(\lambda(\uw))| Q_{\lambda(\uw)}.
\]
Le nombre de $w' \in \WHMj$ tel que $\lambda(\uw) \sim_{\eta_j} w'$ est simplement le cardinal de la classe de $\eta_j$-conjugaison de $\lambda(\uw)$, soit $|\WHMj|/|\Cent_{\eta_j}(\lambda(\uw))|$. On a donc
\begin{align*}
k(w_1w)(\sigma v)  &= \sum_{j=1}^r
\frac{|\Cent_{\gamma_1}(w)|}{|\WHMj|}
\sum_{\unu \in \cN_j} \frac{1}{|W(\unu)|}
\sum_{\substack{\uw \in W(\unu)\\
\kappa(\uw) \sim_{\gamma_1} w}}
\frac{|\WHMj|}{|\Cent_{\eta_j}(\lambda(\uw))|}
 |\Cent_{\eta_j}(\lambda(\uw))| Q_{\lambda(\uw)} \\
&= \sum_{j=1}^r
|\Cent_{\gamma_1}(w)| 
\sum_{\unu \in \cN_j} \frac{1}{|W(\unu)|}
\sum_{\substack{\uw \in W(\unu)\\
\kappa(\uw) \sim_{\gamma_1} w}} Q_{\lambda(\uw)},
\end{align*}
ce qui est la formule cherch\'ee.
\end{proof}

Enfin, nous pouvons combiner les deux lemmes pr\'ec\'edents avec les d\'eveloppements de la Section~\ref{sect:rappels}.

\begin{thm}\label{thm:princ}
Pour tout $E \in \Irr(\WGSE)_\ext$ et tout $g = \sigma v \in G$, on a
\[
\chi_{A_E}(\sigma v) = \sum_{j=1}^r \sum_{\unu \in \cN_j}\, \sum_{E' \in
\Irr(\WHMj)_\ext} \langle E,E'\rangle_\unu\, Q_{E'}(v).
\]
\end{thm}
\begin{proof}
D'apr\`es~\eqref{eqn:AK}, on a
\begin{align*}
\chi_{A_E}(\sigma v)
&= \frac{1}{|\WGSE|} \sum_{w \in \WGSE} \Tr(\gamma_1 w, \tilde E) \chi_{K^{w_1w},\varphi^{w_1w}} \\
&= \frac{1}{|\WGSE|} \sum_{w \in \WGSE} \Tr(\gamma_1 w, \tilde E)
|\Cent_{\gamma_1}(w)| \sum_{j=1}^r
\sum_{\unu \in \cN_j} \frac{1}{|W(\unu)|}
\sum_{\substack{\uw \in W(\unu)\\ \kappa(\uw) \sim_{\gamma_1} w}} Q_{\lambda(\uw)} \\
&= \sum_{j=1}^r \sum_{\unu \in \cN_j} \frac{1}{|W(\unu)|}
\sum_{\uw \in W(\unu)}\ 
\left(\sum_{\substack{w \in \WGSE\\ \kappa(\uw) \sim_{\gamma_1} w}}
\frac{|\Cent_{\gamma_1}(\kappa(\uw))|}{|\WGSE|} \right) 
\Tr(\gamma_1 \kappa(\uw), \tilde E) Q_{\lambda(\uw)}.
\end{align*}
Ici, on a utilis\'e deux cons\'equences du fait que $\kappa(\uw)$ et $w$ sont $\gamma_1$-conjugu\'es: d'une part, $|\Cent_{\gamma_1}(w)| = |\Cent_{\gamma_1}(\kappa(\uw))|$, et d'autre part, on en d\'eduit que les \'el\'ements $\gamma_1\kappa(\uw)$ et $\gamma_1w$ de $\tWGSE$ sont conjugu\'es, et donc que $\Tr(\gamma_1 w, \tilde E) = \Tr(\gamma_1 \kappa(\uw), \tilde E)$.

Ensuite, on peut simplement enlever l'expression entre parenth\`eses ci-dessus: le nombre de $w \in \WGSE$ qui sont $\gamma_1$-conjugu\'es \`a $\kappa(\uw)$ n'est autre que le cardinal de la classe de $\gamma_1$-conjugaison de $\kappa(\uw)$, soit $|\WGSE|/|\Cent_{\gamma_1}(\kappa(\uw))|$. La formule devient donc:
\begin{multline*}
\chi_{A_E}(\sigma v)\\
\begin{aligned}
&= \sum_{j=1}^r \sum_{\unu \in \cN_j} \frac{1}{|W(\unu)|} \sum_{\uw \in W(\unu)}
 \Tr(\gamma_1 \kappa(\uw), \tilde E) Q_{\lambda(\uw)}(v) \\
&= \sum_{j=1}^r \sum_{\unu \in \cN_j} \frac{1}{|W(\unu)|} \sum_{\uw \in W(\unu)}
 \Tr(\gamma_1 \kappa(\uw), \tilde E)
 \sum_{E' \in \Irr(\WHMj)_\ext} \Tr(\gamma' \lambda(\uw)), \tilde E')
 Q_{E'}(v) \\
&= \sum_{j=1}^r \sum_{\unu \in \cN_j}\ \sum_{E' \in \Irr(\WHM)_\ext}
\langle E, E' \rangle_\unu Q_{E'}(v).
\end{aligned}
\end{multline*}
\end{proof}

\section{Classification des Sous-groupes de L\'evi admettant un Faisceau
Caract\`ere Cuspidal}
\label{sect:levi}

Supposons $\bG$ semi-simple et quasi-simple. Nous d\'eterminons dans cette section la liste des sous-groupes de L\'evi de $\bG$ (\`a conjugaison pr\`es) qui admettent un faisceau caract\`ere cuspidal. Il n'y a rien d'original ici: Lusztig a donn\'e tr\`es explicitement en~\cite{LCS} la liste des groupes quasi-simples qui admettent un faisceau caract\`ere cuspidal, et il a \'egalement indiqu\'e en~{\it loc. cit.} un proc\'ed\'e pour d\'eterminer si un groupe r\'eductif donn\'e en admet un ou non. Nous ne faisons qu'effectuer ce proc\'ed\'e. 

Rappelons les \'etapes de ce proc\'ed\'e maintenant. Soit $\bK$ un groupe
r\'eductif, et soit $\chi$ un caract\`ere du groupe des composantes de son
centre $\Cent(\bK)/\Cento(\bK)$. Soit $\hat \bK_\chi^\cusp$ l'ensemble des
faisceaux caract\`eres cuspidaux sur $\bK$, \`a isomorphisme pr\`es, sur
lesquels $\Cent(\bK)/\Cento(\bK)$ agit par $\chi$.
\begin{itemize}
\item Si $\bK$ est semi-simple et quasi-simple, Lusztig a d\'etermin\'e explicitement pour chaque $\chi$ si $\hat \bK_\chi^\cusp$ est vide ou non vide.
\item Si $\bK$ est produit direct des groupes semi-simples et quasi-simples, $\bK = \bK_1 \times \cdots \times \bK_n$, alors $\hat \bK_\chi^\cusp$ est non vide si et seulement si $(\hat\bK_i)_{\chi|_{\bK_i}}^\cusp$ est non vide pour tout $i$.
\item Si $\bK$ est un quotient central d'un produit direct $\tilde \bK$ des groupes semi-simples et quasi-simples, soit $\pi\colon\tilde\bK \to \bK$ l'application quotient. Alors $\hat\bK_\chi^\cusp$ est non vide si et seulement si $\Hat{\Tilde\bK}_{\chi \circ \pi}^\cusp$ l'est.
\item Si $\bK$ est non semi-simple, alors $\hat\bK_\chi^\cusp$ est non vide si et seulement si $(\widehat{\bK/\Cento(\bK)})_\chi^\cusp$ l'est.
\end{itemize}
(En fait, il est \'egalement possible d'obtenir un param\'etrage explicite de $\hat\bK_\chi^\cusp$ au moyen de ce proc\'ed\'e, mais nous n'en aurons pas besoin dans la suite.)

\begin{table}
\rotatebox{90}{\hbox{%
$
\begin{array}{l|l|l}
\bG & \bL & \bM \\
\hline
\SL_{n+1}/\mu_d,\ d \mid (n+1)_{p'} & (\AA_r)^{(n+1)/(r+1)},\
(r+1)d|(n+1)_{p'},\
& (\AA_r)^{(n+1)/(r+1)} \\
\hline
\Spin_{2n+1}\ (p \ne 2) & \BB_{r+s} \times (\AA_1)^{(n-(r+s))/2},\ 2r+1,2s
\in
\triangle & \BB_r \times \DD_s \times (\AA_1)^{(n-(r+s))/2} \\
\SO_{2n+1} & \BB_{r+s},\ 2r+1,2s \in \square\ (p \ne 2) 
& \BB_r \times \DD_s\ (p \ne 2), \\
& \BB_r,\ r \in 2\triangle\ (p = 2) & \BB_r \\
\hline
\Sp_{2n}\ (p \ne 2) & \CC_{r+s},\ r+s \notin 2\Z,\ r,s \in \triangle &
\CC_r \times \CC_s \\
\PSp_{2n} & \CC_{r+s},\ r+s \in 2\Z,\ r,s \in \triangle\ (p \ne 2) & \CC_r
\times \CC_s \\
& \CC_r,\ r \in 2\triangle\ (p = 2) & \CC_r \\
\hline
\Spin_{2n}\ (p \ne 2) & \DD_{r+s} \times (\AA_1)^{(n-(r+s))/2},\ 2r, 2s \in
\triangle & \DD_r \times \DD_s \times (\AA_1)^k \\
\hSpin_{2n},\ n \in 2\Z\ (p \ne 2) & \DD_{r+s} \times
(\AA_1)^{(n-(r+s))/2},\
2r, 2s \in \triangle & \DD_r \times \DD_s \times (\AA_1)^{(n-(r+s))/2}\\
\SO_{2n}\ (p \ne 2) & \DD_{r+s},\ r+s \in 4\Z + 2,\ 2r, 2s \in \square &
\DD_r
\times \DD_s \\
\PSO_{2n} & \DD_{r+s},\ r+s \in 4\Z,\ 2r, 2s \in \square\ (p \ne 2) & \DD_r
\times \DD_s \\
& \DD_r,\ r \in 4\square\ (p = 2) & \DD_r\\
\hline
\EE_6^{\simc}\ (p \ne 3) & (\AA_2)^2 & (\AA_2)^2 \\
          & \EE_6 & \AA_5 \times \AA_1\ (p \ne 2), \EE_6 \\
\EE_6^{\ad} & \DD_4 & (\AA_1)^4\ (p \ne 2), \DD_4\ (p = 2) \\
          & \EE_6 & (\AA_2)^3\ (p \ne 3), \EE_6\ (p = 3) \\
\hline
\EE_7^{\simc}\ (p \ne 2) & (\AA_1)^3 \quad
\text{\small (voir l'explication dans le texte)} &
(\AA_1)^3 \\
        & \EE_7 & \AA_5 \times \AA_2\ (p \ne 3), \EE_7\ (p = 3) \\
\EE_7^{\ad} & \DD_4 & (\AA_1)^4\ (p \ne 2), \DD_4\ (p = 2) \\
        & \EE_6 & (\AA_2)^3\ (p \ne 3), \EE_6\ (p = 3) \\
        & \EE_7 & (\AA_3)^2 \times \AA_1\ (p \ne 2), \EE_7\ (p = 2) \\
\hline
\EE_8 & \DD_4 & (\AA_1)^4\ (p \ne 2), \DD_4\ (p = 2) \\
    & \EE_6 & (\AA_2)^3\ (p \ne 3), \EE_6\ (p = 3) \\
    & \EE_7 & (\AA_3)^2 \times \AA_1\ (p \ne 2), \EE_7\ (p = 2) \\
    & \EE_8\ (p \ne 2) & (\AA_4)^2, \AA_5 \times \AA_2 \times \AA_1, \DD_5 \times
\AA_3, \DD_8, 
    \EE_6 \times \AA_2, \EE_7 \times \AA_1, \EE_8 \\
\hline
\FF_4 & \BB_2 & (\AA_1)^2\ (p \ne 2), \BB_2\ (p = 2) \\
    & \FF_4\ (p \ne 2) & \CC_3 \times \AA_1, \AA_2 \times \AA_2, \AA_3 \times \AA_1,
\BB_4, \FF_4 \\
\hline
\GG_2 & \GG_2 & \AA_1 \times \tilde \AA_1, \AA_2, \GG_2
\end{array}$%
}}
\caption{Sous-groupes de L\'evi admettant un faisceau
caract\`ere cuspidal}\label{tbl:levi}
\end{table}

Les observations suivantes nous seront utiles:

\begin{lem}\label{lem:non-A}
Tout sous-groupe de L\'evi d'un groupe alg\'ebrique quasi-simple poss\`ede au
plus un facteur quasi-simple de type diff\'erent de $A$.
\end{lem}
\begin{proof}
Le graphe de Dynkin d'un groupe quasi-simple et non de type $A$ doit
contenir soit une ar\^ete de multiplicit\'e $2$ ou $3$, soit un n\oe ud de
valence $3$. Chaque graphe de Dynkin simple contient au plus une telle
ar\^ete ou un tel n\oe ud.
\end{proof}

\begin{cor}\label{cor:adjoint-cusp}
Si $\bG$ est quasi-simple et \`a centre connexe, alors tout sous-groupe de
L\'evi $\bL$ admettant un faisceau caract\`ere cuspidal est quasi-simple et
non de type $\AA$.
\end{cor}

Les r\'esultats de la classification sont r\'esum\'es dans la Table~\ref{tbl:levi}.  Pour chaque sous-groupe de L\'evi $\bL$ qui poss\`ede une paire cuspidale $(\Sigma,\cE)$, nous indiquons dans la troisi\`eme colonne le type du groupe $\bM = \Cento_\bL(\sigma)$, o\`u $\sigma \in \Sigma_\semi$.  Lorsque $\bL$ poss\`ede plusieurs paires cuspidales, il y a plusieurs possibilit\'es pour $\bM$.  $\triangle$ d\'esigne l'ensemble des nombres triangulaires, et $\square$ d\'esigne l'ensemble des nombres carr\'es.

Les cas o\`u la caract\'eristique est $2$ et $\bL = \bG$ est de type $\FF_4$ ou $\EE_8$ ne sont pas trait\'es dans la table. Ce sont les seuls cas pour lesquels 
l'hypoth\`ese de nettet\'e (``clean'' au sens de Lusztig) n'est pas encore connue (voir~\cite{ost}).

\subsection{Quotients du groupe lin\'eaire sp\'ecial}

Le centre de $\SL_{n+1}$ est cyclique de cardinal $(n+1)_{p'}$, o\`u
$(n+1)_{p'}$ est le plus grand diviseur de $n+1$ que $p$ ne divise pas.
Soit $d$ un entier qui divise $(n+1)_{p'}$, et notons $\mu_d$ le
sous-groupe cyclique central de cardinal $d$. Tout groupe semi-simple et
quasi-simple de type $\AA_n$ est isomorphe \`a $\SL_{n+1}/\mu_d$ pour un
certain $d$. Ensuite, tout sous-groupe de L\'evi de $\SL_{n+1}$ est de la
forme
\[
S(\GL_{n_1} \times \cdots \times \GL_{n_j})
\qquad\text{o\`u $n_1 + \cdots + n_j = n$,}
\]
et o\`u $S(\cdot)$ signifie le sous-groupe des \'el\'ements \`a
d\'eterminant $1$. Son centre a $\pgcd(n_1, \ldots, n_j)_{p'}$ composantes,
et son image dans $\SL_{n+1}/\mu_d$ a $\pgcd(n_1, \ldots, n_j,
(n+1)_{p'}/d)$ composantes.

Un groupe de type $\AA_n$ admet un faisceau caract\`ere cuspidal si et
seulement si son centre admet un caract\`ere d'ordre $n+1$. Pour que
l'image de $S(\GL_{n_1} \times \cdots \times \GL_{n_j})$ dans
$\SL_{n+1}/\mu_d$ admette un faisceau caract\`ere cuspidal, alors, il faut
et il suffit que $n_i$ divise $\pgcd(n_1, \ldots, n_j, (n+1)_{p'}/d)$ pour
tout $i$. Mais cela implique que
\[
n_1 = \cdots = n_j = \pgcd(n_1, \ldots, n_j, (n+1)_{p'}/d).
\]
Posons $r = n_1 - 1 = \cdots = n_j - 1$. Alors on voit que $r+1$ divise
$(n+1)_{p'}/d$, et que $j = (n+1)/(r+1)$. On conclut qu'un sous-groupe de
L\'evi admet un faisceau caract\`ere cuspidal si et seulement s'il est de
type
\[
\underbrace{\AA_r \times \cdots \AA_r\rlap{,}}_{\text{$(n+1)/(r+1)$
facteurs}}
\qquad\text{o\`u $(r+1) \mid (n+1)_{p'}/d$.}
\]

\subsection{Les groupes classiques}
\label{subsect:cusp-classiques}

En caract\'eristique $2$, tout groupe classique est isomorphe au groupe
adjoint du m\^eme type. Ces groupes-l\`a seront consid\'er\'es dans la
prochaine section; pour le moment, supposons que $p \ne 2$.

Consid\'erons d'abord les groupes sp\'eciaux orthogonaux impairs
$\SO_{2n+1}$. Il est bien connu que tout sous-groupe de L\'evi $\bL$ de
$\SO_{2n+1}$ est de la forme $\SO_{2k+1} \times \GL_{n_1} \times \cdots
\times \GL_{n_j}$, o\`u $(2k+1) + 2n_1 + \cdots + 2n_j = 2n+1$. Mais
$\GL_{n_i}$ n'admet pas de faisceau caract\`ere cuspidal sauf si $n_i = 1$,
et donc pour que $\bL$ en admette un, il doit \^etre de la forme
$\SO_{2k+1} \times \bS$, o\`u $\bS$ est un tore. En particulier, on a
$\bL/\Cento(\bL) \simeq \SO_{2k+1}$, et Lusztig a d\'ecrit en~\cite[\S
23.2(c)]{LCS} des conditions n\'ecessaires sur $k$ pour que $\SO_{2k+1}$
admette un faisceau caract\`ere cuspidal.

Pour $\bG = \Sp_{2n}$, le m\^eme argument permet de se ramener
\`a~\cite[\S 23.2(b)]{LCS}; et pour $\bG = \SO_{2n}$,
\`a~\cite[\S23.2(d)]{LCS}.

\subsection{Les groupes adjoints de type classique}

Selon le Corollaire~\ref{cor:adjoint-cusp}, il suffit de consid\'erer les
sous-groupes de L\'evi $\bL$ qui sont quasi-simples et du m\^eme type que
$\bG$. Si $\bG = \PSp_{m}$ (resp.~$\bG = \PSO_m$), alors il s'ensuit que
$\bL/\Cento(\bL) \simeq \PSp_k$ (resp.~$\bL/\Cento(\bL) \simeq \PSO_k$)
pour un certain $k \le m$ (rappelons que $\bL$ est \`a centre connexe). Si
$p \ne 2$, on se r\'ef\`ere \`a~\cite[\S 23.2(a),(c)]{LCS} pour trouver les
conditions necessaires sur $k$. Si $p = 2$, les r\'esultats analogues se
trouvent en~\cite[\S 22]{LCS}.

\subsection{Les groupes $\Spin$ et $\hSpin$}

Consid\'erons d'abord le groupe $\Spin_{2n+1}$. Son centre est de
cardinal $2$, et donc {\it a priori} il est possible qu'un sous-groupe de
L\'evi contenant des facteurs de type $\AA_1$ puisse admettre un faisceau
caract\`ere cuspidal. Notons $P$ son r\'eseau des poids, et $Q$ son
r\'eseau radiciel. On peut identifier $Q$ avec $\Z^n$,
et $P$ avec le r\'eseau engendr\'e par $Q$ et l'\'el\'ement
\[
\lambda = \textstyle(\frac{1}{2}, \frac{1}{2}, \ldots, \frac{1}{2}).
\]
Posons
\[
e_i = (0, \ldots, 0,1,0,\ldots,0)
\qquad\text{($1$ dans la $i$-\`eme coordonn\'ee)},
\]
et prenons $\{e_1-e_2, e_2-e_3, \ldots, e_{n-1}-e_n, e_n\}$ comme
l'ensemble des racines simples. Soit $\bL$ un sous-groupe de type $\BB_k
\times (\AA_1)^j$. On peut supposer que l'ensemble des racines simples de
$\bL$ est
\begin{multline*}
\{e_1 - e_2\} \cup \{e_3-e_4\} \cup \cdots \cup \{e_{2j-1}-e_{2j}\}
\cup{}\\
\{e_{n-k+1} - e_{n-k+2}, e_{n-k+2}-e_{n-k+3}, \ldots, e_{n-1}-e_n,e_n\},
\end{multline*}
o\`u $n-k+1 > 2j$. Le groupe $\Cent(\bL)/\Cento(\bL)$ poss\`ede un
caract\`ere non trivial si et seulement si un multiple de $\lambda$
appartient au r\'eseau radiciel de $\bL$. Il est donc clair que
$\Cent(\bL)$ est non connexe si et seulement si $2j = n-k$. Ensuite, si
$2j = n-k$, alors $\bL$ admet un faisceau caract\`ere cuspidal si et
seulement si le groupe $\Spin_{2k+1}$ en admet un. Pour ce dernier,
Lusztig a donn\'e les conditions sur $k$ en~\cite[\S 23.2(e)]{LCS}.

Les arguments pour les groupes $\Spin_{2n}$ et $\hSpin_{2n}$ sont
semblables. Pour ceux-ci, on peut identifier le r\'eseau radiciel $Q$ avec
l'ensemble $\{(m_1, \ldots, m_n) \in \Z^n : \sum m_i \in 2\Z\}$. Le
r\'eseau des poids $P$ de $\Spin_{2n}$ est engendr\'e par $Q$ et les deux
\'el\'ements
\[
\lambda = \textstyle(\frac{1}{2}, \frac{1}{2}, \ldots, \frac{1}{2})
\qquad\text{et}\qquad \mu = (0, \ldots, 0, 1).
\]
Le r\'eseau des poids de $\hSpin_{2n}$ est engendr\'e par $Q$ et $\lambda$
seul. L'ensemble des racines simples est $\{e_1-e_2, e_2-e_3, \ldots,
e_{n-1}-e_n, e_{n-1}+e_n\}$. Consid\'erons d'abord un sous-groupe de L\'evi
$\bL$ de type $\DD_k \times (\AA_1)^j$, dont les racines simples sont
\begin{multline*}
\{e_1 - e_2\} \cup \{e_3-e_4\} \cup \cdots \cup \{e_{2j-1}-e_{2j}\}
\cup{}\\
\{e_{n-k+1} - e_{n-k+2}, e_{n-k+2}-e_{n-k+3}, \ldots, e_{n-1}-e_n,e_n\},
\end{multline*}
o\`u $n-k+1 > 2j$. Le poids $\mu$ ne joue aucun r\^ole dans la question,
car le caract\`ere du centre de $\bG$ correspondant est de restriction
nulle aux facteurs de type $\AA_1$. Quant \`a $\lambda$, le m\^eme calcul
qu'on a fait pour $\Spin_{2n+1}$ montre qu'il donne lieu \`a un
caract\`ere non trivial de $\Cent(\bL)/\Cento(\bL)$ si et seulement si $2j
= n - k$. Si $\bG = \Spin_{2n}$ (resp.~$\hSpin_{2n}$), son sous-groupe de
type $\DD_k \times (\AA_1)^{(n-k)/2}$ admet un faisceau caract\`ere
cuspidal si et seulement si $\Spin_{2k}$ (resp.~$\hSpin_{2k}$) en admet
un. Pour ce dernier, voir~\cite[\S 23.2(e),(f)]{LCS}.

Enfin, si $n$ est impair, le centre de $\Spin_{2n}$ est cyclique de
cardinal $4$; en effet, $P/Q$ est engendr\'e par l'image de $\lambda$, et
on a $2\lambda \equiv \mu \pmod Q$. On est donc oblig\'e de consid\'erer
aussi les sous-groupes de L\'evi $\bL$ contenant des facteurs de type
$\AA_3$. Pourtant, le calcul du paragraphe pr\'ec\'edent montre que le
caract\`ere du centre correspondant \`a $\mu$ est toujours de restriction
nulle aux facteurs de type $\AA$. Par cons\'equent, bien que le caract\`ere
correspondant \`a $\lambda$ soit d'ordre $4$, sa restriction aux facteurs
de type $\AA$ n'est que d'ordre $2$. Un tel $\bL$ n'admet donc pas de
faisceau caract\`ere cuspidal.

\subsection{Les groupes adjoints de type exceptionnel}

D'apr\`es le Corollaire~\ref{cor:adjoint-cusp}, il suffit de consid\'erer
les sous-groupes de L\'evi quasi-simples et non de type $\AA$. Tous les tels
groupes sont indiqu\'es dans la Table~\ref{tbl:exc-adj-levi}.

\begin{table}
\[
\begin{array}{c|l}
\bG & \bL \\
\hline
\EE_6 & \DD_4, \DD_5, \EE_6 \\
\EE_7 & \DD_4, \DD_5, \DD_6, \EE_6, \EE_7 \\
\EE_8 & \DD_4, \DD_5, \DD_6, \DD_7, \EE_6, \EE_7, \EE_8 \\
\FF_4 & \BB_2, \BB_3, \CC_3, \FF_4 \\
\GG_2 & \GG_2
\end{array}
\]
\caption{Sous-groupes de L\'evi quasi-simples et non de type $\AA$ dans les
groupes exceptionnels}\label{tbl:exc-adj-levi}
\end{table}

Parmi ceux-ci, les groupes adjoints de type $\BB_3$, $\CC_3$, $\DD_5$,
$\DD_6$, $\DD_7$ n'admettent pas de faisceaux caract\`eres
cuspidaux~\cite[\S 22,\S 23.2(a),(c)]{LCS}. Tous les autres en admettent au
moins un. Les groupes de type $\BB_2$ et $\DD_4$ en caract\'eristique $2$ sont trait\'es en~\cite[\S 22]{LCS}, et en caract\'eristique impaire dans la Proposition~23.2(c) de~{\it loc. cit.}.  Pour la liste des $\bM$ possibles dans $\EE_6$, $\EE_7$, $\EE_8$, $\FF_4$, $\GG_2$, respectivement, voir les
Propositions~20.3(a), 20.3(c), 21.2, 21.3, 20.6
de~\cite{LCS}.

\subsection{Le groupe simplement connexe de type $\EE_6$}

Si $p = 3$, $\EE_6^\simc$ est isomorphe \`a $\EE_6^\ad$. Supposons donc
que $p \ne 3$. Puisque le centre de $\EE_6^\simc$ est de cardinal $3$, on sait que pour tout sous-groupe de L\'evi $\bL$, le cardinal de $\Cent(\bL)/\Cento(\bL)$ \'egale soit $1$, soit $3$.  Parmi les sous-groupes figurant dans la
Table~\ref{tbl:exc-adj-levi}, ceux de type $\DD_4$ et $\DD_5$ sont \`a
centre connexe (parce qu'un groupe quasi-simple de type $\DD$ ne peut pas avoir un centre \`a $3$ composantes) et donc ont d\'ej\`a \'et\'e trait\'e. Pour $\EE_6^\simc$ lui-m\^eme, voir~\cite[Proposition~20.3]{LCS}.

Nous devons maintenant consid\'erer les sous-groupes de L\'evi non quasi-simples contenant un facteur de type $\AA_2$. Il y en a deux, de types $\AA_2$ et $\AA_2 \times \AA_2$. Un calcul semblable \`a ceux que l'on a fait pour les groupes $\Spin$ montre que le sous-groupe de L\'evi de type $\AA_2$ est \`a centre connexe. (Il suffit de v\'erifier, d'apr\`es les descriptions en~\cite{Bour} des r\'eseaux de poids et des r\'eseaux radiciels, qu'il n'existe pas de poids qui n'est pas dans le r\'eseau radiciel de $\AA_2$ mais dont un multiple y est). Par contre, le centre de $\AA_2 \times \AA_2$ est \`a $3$ composantes, et donc ce groupe-ci admet un faisceau caract\`ere cuspidal.

\subsection{Le groupe simplement connexe de type $\EE_7$}

Si $p = 2$, $\EE_7^\simc$ est isomorphe \`a $\EE_7^\ad$. Supposons donc
que $p \ne 2$. Le centre de $\EE_7^\simc$ est alors de cardinal $2$. Il n'est pas donc aussi facile qu'en $\EE_6^\simc$ de conclure que les divers sous-groupes propres de L\'evi figurant dans la Table~\ref{tbl:exc-adj-levi} sont \`a centre connexe. Il faut plut\^ot v\'erifier ce fait dans chaque cas par un calcul dans le r\'eseau des poids, en utilisant les donn\'ees de~\cite{Bour}.

Des calculs semblables montrent que le seul sous-groupe de L\'evi (\`a conjugaison pr\`es) contenant des facteurs de type $\AA_1$ et poss\'edant une paire cuspidale est celui de type $\AA_1 \times \AA_1 \times \AA_1$ qui correspond au graphe suivant:
\[
\vcenter{\hbox{\small $\xymatrix@=4pt{
*{}&*{}& *{\bullet} \ar@{-}[d] \\
*{\circ}\ar@{-}[r] &
*{\circ}\ar@{-}[r] & *{\circ}\ar@{-}[r] &
*{\bullet}\ar@{-}[r] & *{\circ}\ar@{-}[r] & *{\bullet}}$}}
\]
($\EE_7^\simc$ poss\`ede plusieurs sous-groupes de L\'evi non conjugu\'es de
type $\AA_1 \times \AA_1 \times \AA_1$. Les autres sous-groupes de ce type
n'ont pas de faisceaux caract\`eres cuspidaux). 

\section{Classes de Conjugaison de Sous-groupes de L\'evi dans le
Centralisateur d'un \'El\'ement semi-simple}
\label{sect:sousgp}

On garde les notations de la section pr\'ec\'edente: $\bG$ est un groupe r\'eductif, $\bL = \Cent_\bG(\bT_1)$ est un sous-groupe de L\'evi qui poss\`ede une paire cuspidale $(\Sigma,\cE)$, et $\cM = \Cento_\bL(\sigma)$, o\`u $\sigma \in \Sigma_\semi$. Rappelons la d\'efinition de l'ensemble $\cM$ qu'on a introduit \`a la Section~\ref{sect:notations}:
\[
\cM = \{ \bM_x : x \in \bG,\ x^{-1} \sigma x \in \Sigma_\semi \}
\qquad\text{o\`u}\qquad
\bM_x = \Cento_\bH(x\bT_1 x^{-1}).
\]
\`A ce moment-l\`a, nous avons remarqu\'e que $\cM$ se r\'epartit en un nombre fini de classes de conjugaison sous $\bH$.

\begin{remark}
\`A la Section~\ref{sect:notations}, la notation $\sigma$ d\'esignait un \'el\'ement de $\bG$ conjugu\'e \`a un \'el\'ement de $\Sigma_\semi$. Ici, on l'a suppos\'e dans $\Sigma_\semi$. Il est clair qu'aucune perte de g\'en\'eralit\'e n'en r\'esulte.
\end{remark}

Cette section est consacr\'ee \`a la preuve du th\'eor\`eme suivant:

\begin{thm}\label{thm:m-conj}
Si $\bG$ est semi-simple, quasi-simple, et diff\'erent de $\PSp_{2n}$, $\PSO_{2n}$, $\hSpin_{2n}$ et $\EE_7^\simc$, alors tous les membres de $\cM$ sont conjugu\'es sous $\bH$. Si $\bG$ est l'un de ces quatres groupes, alors $\cM$ se r\'epartit en une ou deux classes de conjugaison sous $\bH$.
\end{thm}

Il est \`a rappeler (voir la Remarque~\ref{rmk:m-unique}) que tous les membres de $\cM$ sont conjugu\'es sous $\bG$ et donc isomorphes.

\begin{remark}\label{rmk:m-conj-hyp}
Soit $x \in \bG$ est tel que $x^{-1} \sigma x \in \Sigma_\semi$, et posons $\sigma' = x^{-1} \sigma x$. On sait, d'apr\`es la Remarque~\ref{rmk:m-unique}, que $\sigma$ et $\sigma'$ sont conjugu\'es (dans $\bL$) \`a multiplication par un \'el\'ement de $\bT_1$ pr\`es. Il y a donc un $f \in \bL$ et un $z \in \bT_1$ tel que $f^{-1}\sigma'f = z\sigma$, ou autrement dit, $(xf)^{-1}\sigma (xf) = z \sigma$. Puisque $f$ centralise $\bT_1$, il est clair $xf \bT_1 (xf)^{-1} = x \bT_1 x^{-1}$, et donc que $\bM_{xf} = \bM_x$. On ne s'int\'eresse qu'\`a la classe de conjugaison de ce dernier groupe, et donc on peut supposer, sans perte de g\'en\'eralit\'e, que $x^{-1} \sigma x = z \sigma$, avec $z \in \bT_1$.
\end{remark}

\subsection{Les cas triviaux}

Si $\bL = \bG$, alors on a $\bH = \bM = \bM_x$ pour tout $x$, et il n'y a rien \`a d\'emontrer. En particulier, le th\'eor\`eme est donc vrai dans les cas suivants:
\[
\begin{array}{c|c}
\bG = \bL & \bM \\
\hline
\EE_6 & \AA_5 \times \AA_1, (\AA_2)^3, \EE_6 \\
\EE_7 & \AA_5 \times \AA_2, (\AA_3)^2 \times \AA_1, \EE_7 \\
\EE_8 & (\AA_4)^2, \AA_5 \times \AA_2 \times \AA_1, \DD_5 \times
\AA_3, \DD_8, \EE_6 \times \AA_2, \EE_7 \times \AA_1, \EE_8 \\
\FF_4 & \CC_3 \times \AA_1, \AA_2 \times \AA_2, \AA_3 \times \AA_1, \BB_4, \FF_4 \\
\GG_2 & \AA_1 \times \AA_1, \AA_2, \GG_2
\end{array}
\]

\subsection{Graphes de Dynkin}
\label{subsect:m-graphe}

Soient $\Delta_\bG$, $\Delta_\bL$, $\Delta_\bH$, $\Delta_\bM$ et $\Delta_{\bM_x}$ les graphes de Dynkin des groupes correspondants, et soient $\tilde\Delta_\bG$ et $\tilde\Delta_\bL$ les graphes de Dynkin compl\'et\'es de $\bG$ et de $\bL$. Rappelons que le graphe de Dynkin d'un sous-groupe de L\'evi (resp.~du centralisateur 
d'un \'el\'ement semi-simple) peut \^etre identifi\'e \`a un sous-graphe, unique 
\`a conjugaison sous le groupe de Weyl pr\`es, du graphe de Dynkin (resp.~graphe de Dynkin compl\'et\'e) du groupe de d\'epart. On a donc les inclusions suivantes:
\[
\Delta_\bL \subset \Delta_\bG
\qquad\text{et}\qquad
\Delta_{\bM_x} \subset \Delta_\bH \subset \tilde\Delta_\bG,
\]
On sait d\'ej\`a que $\Delta_\bM \simeq \Delta_{\bM_x}$ pour tout $x$. Pour montrer que tous les $\bM_x$ sont conjugu\'es dans $\bH$, il suffit de montrer que les sous-graphes $\Delta_{\bM_x}$ de $\Delta_\bH$ sont conjugu\'es par le groupe de Weyl de $\bH$. En particulier, si tous les sous-graphes de $\Delta_\bH$ isomorphes \`a $\Delta_\bM$ soient conjugu\'es, alors le th\'eor\`eme s'ensuit. 

Il arrive parfois que $\tilde\Delta_\bG$ (et donc $\Delta_\bH$) ne contienne 
qu'un seul sous-graphe isomorphe \`a $\Delta_\bM$. Dans ces cas-l\`a, il n'y a 
rien \`a d\'emontrer, et le r\'esultat est imm\'ediat.

\subsection{Les groupes de type $\AA$ et les groupes exceptionnels}
\label{subsect:m-a-exc}

Si $\bG$ est de type $\AA_n$, alors tout centralisateur d'un \'el\'ement
semi-simple est en fait de L\'evi, et on peut donc se restreindre \`a
consid\'erer le graphe de Dynkin non compl\'et\'e $\Delta_\bG$. $\bM$ est de
type $(\AA_{r-1})^{(n+1)/r}$. Si l'on note $\alpha_1, \ldots, \alpha_n$ les
n\oe uds de $\Delta_\bG$, il est clair que l'unique sous-graphe de type
$(\AA_{r-1})^{(n+1)/r}$ est celui qui contient les n\oe uds
\[
\alpha_1, \ldots, \alpha_{r-1}; \alpha_{r+1}, \ldots, \alpha_{2r-1}; \ldots; \alpha_{(n+1)/r -r +1}, \ldots, \alpha_n.
\]

Supposons maintenant que $(\bG, \bL, \bM)$ est l'un des triplets qui figurent dans la Table~\ref{tbl:exc-gp}. Dans les cas o\`u il n'y a aucune mention sous l'en-t\^ete ``remarque,'' le graphe $\tilde\Delta_\bG$ ne contient qu'un seul sous-graphe isomorphe \`a $\Delta_\bM$.

Dans les trois cas qui portent la mention $(*)$, $\tilde\Delta_\bG$ contient d'autres sous-graphes isomorphes \`a $\Delta_\bM$, et il faut donc faire un argument suppl\'ementaire. Consid\'erons le cas o\`u $\bG$ et de type $\EE_7$ et $\bM$ de type $(\AA_1)^4$. Le groupe ne peut pas \^etre un sous-groupe de L\'evi de $\bG$: on sait que $\bM$ admet un faisceau caract\`ere cuspidal, mais selon la Table~\ref{tbl:levi}, $\bG$ n'a pas de sous-groupe de L\'evi isog\`ene \`a $\bM$ qui en admet un. Il est facile de faire la liste de tous les sous-graphes de $\tilde\Delta_\bG$ de type $(\AA_1)^4$:
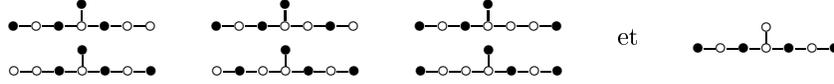
\[
\begin{array}{@{}c@{}}
\vcenter{\hbox{\small $\xymatrix@=4pt{
*{}&*{}&*{}& *{\bullet} \ar@{-}[d] \\
*{\bullet}\ar@{-}[r] & *{\circ}\ar@{-}[r] &
*{\bullet}\ar@{-}[r] & *{\circ}\ar@{-}[r] &
*{\bullet}\ar@{-}[r] & *{\circ}\ar@{-}[r] & *{\circ}}$}}_{\strut}
\\
\vcenter{\hbox{\small $\xymatrix@=4pt{
*{}&*{}&*{}& *{\bullet} \ar@{-}[d] \\
*{\circ}\ar@{-}[r] & *{\circ}\ar@{-}[r] &
*{\bullet}\ar@{-}[r] & *{\circ}\ar@{-}[r] &
*{\bullet}\ar@{-}[r] & *{\circ}\ar@{-}[r] & *{\bullet}}$}}
\end{array}
\qquad
\begin{array}{@{}c@{}}
\vcenter{\hbox{\small $\xymatrix@=4pt{
*{}&*{}&*{}& *{\bullet} \ar@{-}[d] \\
*{\bullet}\ar@{-}[r] & *{\circ}\ar@{-}[r] &
*{\bullet}\ar@{-}[r] & *{\circ}\ar@{-}[r] &
*{\circ}\ar@{-}[r] & *{\bullet}\ar@{-}[r] & *{\circ}}$}}_{\strut}
\\
\vcenter{\hbox{\small $\xymatrix@=4pt{
*{}&*{}&*{}& *{\bullet} \ar@{-}[d] \\
*{\circ}\ar@{-}[r] & *{\bullet}\ar@{-}[r] &
*{\circ}\ar@{-}[r] & *{\circ}\ar@{-}[r] &
*{\bullet}\ar@{-}[r] & *{\circ}\ar@{-}[r] & *{\bullet}}$}}
\end{array}
\qquad
\begin{array}{@{}c@{}}
\vcenter{\hbox{\small $\xymatrix@=4pt{
*{}&*{}&*{}& *{\bullet} \ar@{-}[d] \\
*{\bullet}\ar@{-}[r] & *{\circ}\ar@{-}[r] &
*{\bullet}\ar@{-}[r] & *{\circ}\ar@{-}[r] &
*{\circ}\ar@{-}[r] & *{\circ}\ar@{-}[r] & *{\bullet}}$}}_{\strut}
\\
\vcenter{\hbox{\small $\xymatrix@=4pt{
*{}&*{}&*{}& *{\bullet} \ar@{-}[d] \\
*{\bullet}\ar@{-}[r] & *{\circ}\ar@{-}[r] &
*{\circ}\ar@{-}[r] & *{\circ}\ar@{-}[r] &
*{\bullet}\ar@{-}[r] & *{\circ}\ar@{-}[r] & *{\bullet}}$}}
\end{array}
\qquad\text{et}\qquad
\vcenter{\hbox{\small $\xymatrix@=4pt{
*{}&*{}&*{}& *{\circ} \ar@{-}[d] \\
*{\bullet}\ar@{-}[r] & *{\circ}\ar@{-}[r] &
*{\bullet}\ar@{-}[r] & *{\circ}\ar@{-}[r] &
*{\bullet}\ar@{-}[r] & *{\circ}\ar@{-}[r] & *{\bullet}}$}}
\]
Les six premiers sous-graphes sont tous conjugu\'es sous le groupe de Weyl de $\bG$, et chacun des quatre premiers sous-graphes est contenu dans un sous-graphe de type $\EE_7$, et correspond donc \`a un sous-groupe de L\'evi. Seul le dernier correspond \`a un sous-groupe qui n'est pas de L\'evi; celui-l\`a doit \^etre \'egal \`a $\Delta_\bM$. 

Un argument semblable permet de traiter les autres cas marqu\'es~$(*)$: dans chaque cas, on trouve d'apr\`es la Table~\ref{tbl:levi} que $\bM$ n'est pas un sous-groupe de L\'evi de $\bG$, et que $\tilde\Delta_\bG$ n'a qu'un seul sous-graphe isomorphe \`a $\Delta_\bM$ dont le sous-groupe correspondant n'est pas de L\'evi. Ce sous-graphe-l\`a est donc forc\'ement \'egal et \`a $\Delta_\bM$. En particulier, $\Delta_\bM$ est seul dans sa classe de conjugaison sous le groupe de Weyl $\bG$, et donc sous celui de $\bH$. 


Le th\'eor\`eme est maintenant d\'emontr\'e dans tous les cas figurant dans
la Table~\ref{tbl:exc-gp}.

\begin{table}
\[
\begin{array}{|c|c|c|cc|l|}
\bG & \bL & \bM & \text{graphe} && \text{remarque} \\
\hline
\AA_n & (\AA_{r-1})^{(n+1)/r} & (\AA_{r-1})^{(n+1)/r} &
\multicolumn{2}{l|}{\vcenter{\hbox{\tiny $\xymatrix@=2pt{
*{\bullet}\ar@{.}[rr] &&
*{\bullet}\ar@{-}[r] & *{\circ}\ar@{-}[r] &
*{\bullet}\ar@{.}[rr] &&
*{\bullet}\ar@{-}[r] & *{\circ}\ar@{-}[r] &
*{\ }\ar@{.}[rrr] &&&
*{\ }\ar@{-}[r] & *{\circ}\ar@{-}[r] &
*{\bullet}\ar@{..}[rr] &&
*{\bullet}}$}}}
& \\
\hline
\EE_6 & \DD_4\strut & \AA_1 \times \AA_1 \times \AA_1 \times \AA_1 &
\vcenter{\hbox{\tiny $\xymatrix@=2pt{
*{}&*{}& *{\bullet} \ar@{-}[d] \\
*{}&*{}& *{\circ} \ar@{-}[d] \\
*{\bullet}\ar@{-}[r] & *{\circ}\ar@{-}[r] &
*{\bullet}\ar@{-}[r] & *{\circ}\ar@{-}[r] & *{\bullet}}$}}
& \raisebox{1pt}{\hbox{\strut}}& \\
    & \DD_4 & \DD_4 &
\vcenter{\hbox{\tiny $\xymatrix@=2pt{
*{}&*{}& *{\circ} \ar@{-}[d] \\
*{}&*{}& *{\bullet} \ar@{-}[d] \\
*{\circ}\ar@{-}[r] & *{\bullet}\ar@{-}[r] &
*{\bullet}\ar@{-}[r] & *{\bullet}\ar@{-}[r] & *{\circ}}$}}
& \raisebox{3pt}{\hbox{\strut}} & \\
\hline
\EE_7 & \DD_4 & \AA_1 \times \AA_1 \times \AA_1 \times \AA_1 &
\vcenter{\hbox{\tiny $\xymatrix@=2pt{
*{}&*{}&*{}& *{\circ} \ar@{-}[d] \\
*{\bullet}\ar@{-}[r] & *{\circ}\ar@{-}[r] &
*{\bullet}\ar@{-}[r] & *{\circ}\ar@{-}[r] &
*{\bullet}\ar@{-}[r] & *{\circ}\ar@{-}[r] & *{\bullet}}$}}
&& (*)\\
    & \DD_4 & \DD_4 &
\vcenter{\hbox{\tiny $\xymatrix@=2pt{
*{}&*{}&*{}& *{\bullet} \ar@{-}[d] \\
*{\circ}\ar@{-}[r] & *{\circ}\ar@{-}[r] &
*{\bullet}\ar@{-}[r] & *{\bullet}\ar@{-}[r] &
*{\bullet}\ar@{-}[r] & *{\circ}\ar@{-}[r] & *{\circ}}$}}
&& \\
    & \EE_6 & \AA_2 \times \AA_2 \times \AA_2 &
\vcenter{\hbox{\tiny $\xymatrix@=2pt{
*{}&*{}&*{}& *{\bullet} \ar@{-}[d] \\
*{\bullet}\ar@{-}[r] & *{\bullet}\ar@{-}[r] &
*{\circ}\ar@{-}[r] & *{\bullet}\ar@{-}[r] &
*{\circ}\ar@{-}[r] & *{\bullet}\ar@{-}[r] & *{\bullet}}$}}
&& \\
    & \EE_6 & \EE_6 &
\vcenter{\hbox{\tiny $\xymatrix@=2pt{
*{}&*{}&*{}& *{\bullet} \ar@{-}[d] \\
*{\circ}\ar@{-}[r] & *{\bullet}\ar@{-}[r] &
*{\bullet}\ar@{-}[r] & *{\bullet}\ar@{-}[r] &
*{\bullet}\ar@{-}[r] & *{\bullet}\ar@{-}[r] & *{\circ}}$}}
&& \\
\hline
\EE_8 & \DD_4 & \AA_1 \times \AA_1 \times \AA_1 \times \AA_1 &
\vcenter{\hbox{\tiny $\xymatrix@=2pt{
*{}&*{}& *{\bullet} \ar@{-}[d] \\
*{\circ}\ar@{-}[r] & *{\circ}\ar@{-}[r] &
*{\circ}\ar@{-}[r] & *{\bullet}\ar@{-}[r] & *{\circ}\ar@{-}[r] &
*{\bullet}\ar@{-}[r] & *{\circ}\ar@{-}[r] & *{\bullet}}$}}
&& (*)\\
    & \DD_4 & \DD_4 &
\vcenter{\hbox{\tiny $\xymatrix@=2pt{
*{}&*{}& *{\bullet} \ar@{-}[d] \\
*{\circ}\ar@{-}[r] & *{\bullet}\ar@{-}[r] &
*{\bullet}\ar@{-}[r] & *{\bullet}\ar@{-}[r] & *{\circ}\ar@{-}[r] &
*{\circ}\ar@{-}[r] & *{\circ}\ar@{-}[r] & *{\circ}}$}}
&& \\
    & \EE_6 & \AA_2 \times \AA_2 \times \AA_2 &
\vcenter{\hbox{\tiny $\xymatrix@=2pt{
*{}&*{}& *{\circ} \ar@{-}[d] \\
*{\bullet}\ar@{-}[r] & *{\bullet}\ar@{-}[r] &
*{\circ}\ar@{-}[r] & *{\bullet}\ar@{-}[r] & *{\bullet}\ar@{-}[r] &
*{\circ}\ar@{-}[r] & *{\bullet}\ar@{-}[r] & *{\bullet}}$}}
&& \\
    & \EE_6 & \EE_6 &
\vcenter{\hbox{\tiny $\xymatrix@=2pt{
*{}&*{}& *{\bullet} \ar@{-}[d] \\
*{\bullet}\ar@{-}[r] & *{\bullet}\ar@{-}[r] &
*{\bullet}\ar@{-}[r] & *{\bullet}\ar@{-}[r] & *{\bullet}\ar@{-}[r] &
*{\circ}\ar@{-}[r] & *{\circ}\ar@{-}[r] & *{\circ}}$}}
&& \\
    & \EE_7 & \AA_3 \times \AA_3 \times \AA_1 &
\vcenter{\hbox{\tiny $\xymatrix@=2pt{
*{}&*{}& *{\bullet} \ar@{-}[d] \\
*{\bullet}\ar@{-}[r] & *{\circ}\ar@{-}[r] &
*{\bullet}\ar@{-}[r] & *{\bullet}\ar@{-}[r] & *{\circ}\ar@{-}[r] &
*{\bullet}\ar@{-}[r] & *{\bullet}\ar@{-}[r] & *{\bullet}}$}}
&& \\
    & \EE_7 & \EE_7 &
\vcenter{\hbox{\tiny $\xymatrix@=2pt{
*{}&*{}& *{\bullet} \ar@{-}[d] \\
*{\bullet}\ar@{-}[r] & *{\bullet}\ar@{-}[r] &
*{\bullet}\ar@{-}[r] & *{\bullet}\ar@{-}[r] & *{\bullet}\ar@{-}[r] &
*{\bullet}\ar@{-}[r] & *{\circ}\ar@{-}[r] & *{\circ}}$}}
&& \\
\hline
\FF_4 & \BB_2 & \AA_1 \times \AA_1 &
\vcenter{\hbox{\tiny $\xymatrix@=2pt{
*{\bullet} \ar@{-}[r] & *{\circ} \ar@{-}[r] & *{\bullet}
\ar@2{-}[r]|{\scriptstyle \rangle} \ar@2{-}[r] &
*{\circ} \ar@{-}[r] & *{\circ}}$}}
&& (*) \\
      & \BB_2 & \BB_2 &
\vcenter{\hbox{\tiny $\xymatrix@=2pt{
*{\circ} \ar@{-}[r] & *{\circ} \ar@{-}[r] & *{\bullet}
\ar@2{-}[r]|{\scriptstyle \rangle} \ar@2{-}[r] &
*{\bullet} \ar@{-}[r] & *{\circ}}$}}
&& \\
\hline
\end{array}
\]
\caption{Groupes exceptionnels}\label{tbl:exc-gp}
\end{table}

\subsection{Les groupes classiques}

Supposons maintenant que $\bG$ est l'un des groupes $\SO_{2n+1}$, $\Sp_{2n}$, ou $\SO_{2n}$. Dans les groupes classiques, nous pouvons tirer profit du fait que tout
sous-groupe de L\'evi se d\'ecompose en produit direct de son sous-groupe
d\'eriv\'e et d'un tore. En particulier, on a
\[
\bL \simeq \bL_\der \times \bT_1.
\]
De plus, $\bL$ se plonge dans un sous-groupe r\'eductif
\[
\bL_\der \times \bL_1 \subset \bG
\]
o\`u $\bL_1$ est semi-simple, quasi-simple et du m\^eme type que $\bG$, et o\`u $\bT_1$ est un tore maximal de $\bL_1$. Le groupe $\bL_1$ se d\'ecrit comme suit: si $\bG = \SO_N$ (resp.~$\bG = \Sp_N$), alors il y a un entier positif $M \le N$ tel que $\bL_\der \simeq \SO_M$ (resp.~$\bL_\der \simeq \Sp_M$), et on a que $\bL_1 \simeq \SO_{N - M}$ (resp.~$\bL_1 \simeq \Sp_{N - M}$). Il est \`a noter que $\bL_\der \times \bL_1$ n'est pas en g\'en\'eral de L\'evi.

Si $y$ est un \'el\'ement d'un groupe classique, on note $E(y)$ l'ensemble (avec multiplicit\'es) de ses valeurs propres. Si l'on d\'ecompose un \'el\'ement $y \in \bL_\der \times \bL_1$ en un produit $y = d \times t$, o\`u $d \in \bL_\der$ et $t \in \bL_1$, alors $E(y) = E(d) \cup E(t)$. Enfin, rappelons que dans tout groupe classique quasi-simple, deux \'el\'ements semi-simples sont conjugu\'es si et seulement s'ils ont les m\^emes valeurs propres avec les m\^emes multiplicit\'es.

Imposons l'hypoth\`ese de la Remarque~\ref{rmk:m-conj-hyp}, et reprenons
ses notations: $\sigma' = x^{-1}\sigma x = z\sigma$, o\`u $z \in \bT_1$.
\'Ecrivons $\sigma$ comme un produit $d\cdot t$, o\`u $d \in \bL_\der$ et
$t \in \bT_1$, et de m\^eme pour $\sigma' = d'\cdot t'$. Puisque de telles
d\'ecompositions sont uniques, l'\'egalit\'e $\sigma' = z \sigma$ implique
que
\[
d' = d
\qquad\text{et}\qquad
t' = zt.
\]
Puisque $\sigma$ et $\sigma'$ sont conjugu\'es, on sait que $E(\sigma) = E(\sigma')$, et puis il s'ensuit que $E(t) = E(zt)$. Il existe donc un $k \in \Norm_{\bL_1}(\bT_1)$ tel que $k t k^{-1} = zt$. Puisque $k$ commute avec $\bL_\der$, on a maintenant
\[
k\sigma k^{-1} = k dt k^{-1} = z dt = \sigma'.
\]
Posons $h' = xk$; alors $h' \sigma h'^{-1} = \sigma$. D'autre part, on a que $h' \bT_1 h'^{-1} = x\bT_1 x^{-1}$ (car $k$ normalise $\bT_1$), et donc
\begin{equation}\label{eqn:m-h-classique}
h' \bM h'^{-1} = \bM_x.
\end{equation}
L'\'el\'ement $h'$ appartient au groupe \'eventuellement non connexe $\Cent_\bG(\sigma)$. Il reste de d\'emontrer qu'on peut remplacer $h'$ par un \'el\'ement de $\bH$.

Consid\'erons le groupe non connexe $\Or_m$. Il est clair que pour tout tore $\bS \subset \SO_m$, il existe un \'el\'ement dans la composante non neutre de $\Or_m$ qui commute avec $\bS$. Par cons\'equent, si $\bS$ est un tore dans un produit quelconque des $\Or_m$, des $\Sp_m$, et des $\GL_m$, alors toute composante contient un \'el\'ement commutant avec $\bS$. Le groupe $\Cent_\bG(\sigma)$ est un 
sous-groupe d'un tel produit, et donc dans la composante de $\Cent_\bG(\sigma)$ contenant $h'$, il existe un \'el\'ement $r$ qui centralise le tore $x^{-1}\bT_1 x$. Cette derni\`ere condition implique que $r$ normalise $\bM_x$. Ensuite, posons $h = r^{-1} h'$. Cet \'el\'ement est forc\'ement dans la
composante neutre de $\Cent_\bG(\sigma)$, \ie dans $\bH$. Il s'ensuit maintenant de~\eqref{eqn:m-h-classique} que $h \bM h^{-1} = \bM_x$. $\bM$ et $\bM_x$ sont donc conjugu\'es sous $\bH$.

\begin{remark}
Au cours de cette preuve, l'hypoth\`ese que $\bL$ et $\bM$ admettent des faisceaux caract\`eres cuspidaux n'a jou\'e aucun r\^ole. Cela nous aidera plus tard \`a traiter les groupes simplement connexes de type classique.
\end{remark}

\subsection{Les groupes adjoints de type classique}

Supposons que $\bG$ est l'un des groupes $\PSp_{2n}$ ou $\PSO_{2n}$. La preuve pour ces groupes-ci consiste \`a se ramener au cas des groupes classiques. Posons $\tilde\bG = \Sp_{2n}$ ou $\SO_{2n}$, respectivement, et soit $\pi\colon\tilde\bG \to \bG$ l'application quotient naturelle. Soit $\tilde \bT_1$ la composante neutre de $\pi^{-1}(\bT_1)$, et posons $\tilde\bL = \Cent_{\tilde\bG}(\tilde\bT_1)$ et $\tilde\Sigma_\semi = \pi^{-1}(\Sigma_\semi)$. Ensuite, choisissons un point $\dot\sigma \in \pi^{-1}(\sigma)$, ainsi qu'un $\dot x \in \pi^{-1}(x)$ pour tout $x \in \bG$ tel que $x^{-1}\sigma x \in \Sigma_\semi$. On a donc $\dot x^{-1}\dot \sigma \dot x \in \tilde \Sigma_\semi$. Enfin, posons $\tilde\bM_{\dot x} = \Cent_{\tilde\bH}(\dot x\tilde\bT_1\dot x^{-1})$. L'ensemble $\tilde\Sigma_\semi$ est soit connexe, soit \`a deux composantes. Essayons d'abord de mieux le comprendre.

Soit $C \subset \bL/\bT_1$ la classe de conjugaison dont $\Sigma_\semi$ est l'image r\'eciproque, et consid\'erons l'application $q\colon\tilde\bL/\tilde \bT_1 \to \bL/\bT_1$. L'ensemble $\tilde C = q^{-1}(C)$ contient une ou deux classes de 
conjugaison et une ou deux composantes. D'une part, toute classe de conjugaison dans $\tilde\bL/\tilde \bT_1$ est connexe; d'autre part, la r\'eunion de deux classes de conjugaison de m\^eme dimension est forc\'ement non connexe. On conclut que chaque composante de $\tilde C$ est une classe de conjugaison. Ensuite, l'application $r\colon\tilde\bL \to \tilde\bL/\tilde\bT_1$ \'etant \`a noyau connexe, on voit que l'op\'eration d'image r\'eciproque sous $r$ pr\'eserve le nombre de composantes. Donc chaque composante de $\tilde \Sigma_\semi = r^{-1}(\tilde C)$ est l'image r\'eciproque d'une seule classe de conjugaison dans $\tilde \bL/\tilde \bT_1$.

Si $\tilde\Sigma_\semi$ est connexe, alors tous les $\tilde\bM_{\dot x}$ sont 
conjugu\'es sous $\tilde\bH$ (car le th\'eor\`eme est d\'ej\`a \'etabli pour 
$\tilde\bG$), et donc leurs images $\bM_x$ sont conjugu\'es sous $\bH$.

En revanche, si $\tilde\Sigma_\semi$ est r\'eunion de deux composantes, notons-les $\tilde\Sigma_\semi^+$ et $\tilde\Sigma_\semi^-$. Ensuite, definissons deux sous-ensemble $\cM^+, \cM^- \subset \cM$ par
\[
\cM^\pm = \{ \bM_x : \dot x^{-1}\dot \sigma \dot x \in \tilde\Sigma_\semi^\pm \}.
\]
Il est \`a noter que cette r\'epartition de $\cM$ en deux sous-ensembles est bien d\'efinie, \ie ind\'ependante des choix des \'el\'ements $\dot x$ et $\dot \sigma$. Tous les membres de $\cM^+$ (resp.~$\cM^-$) sont conjugu\'es sous $\bH$, puisque les $\tilde\bM_{\dot x}$ correspondants sont conjugu\'es sous $\tilde\bH$.

Parfois, les membres de $\cM^+$ et de $\cM^-$ deviennent conjugu\'es sous $\bH$, mais il est \'egalement possible qu'ils restent non conjugu\'es. Par exemple, supposons que $\tilde\bG = \SO_{4r}$ et $\tilde \bM = \tilde\bL \simeq \SO_{2r} \times \tilde \bT_1$. Soit $\dot\sigma$ l'\'el\'ement $(-1,1) \in \SO_{2r} \times \tilde\bT_1$. Alors $\tilde\bH \simeq \SO_{2r} \times \SO_{2r}$. Remarquons que les \'el\'ements $\dot\sigma$ et $-\dot\sigma$ ont les m\^emes valeurs propres (avec les m\^emes multiplicit\'es); ils sont donc conjugu\'es dans $\tilde\bG$. En effet, soit $\dot x \in \tilde\bG$ une matrice de permutation telle que la conjugaison par $\dot x$ \'echange les deux facteurs de $\bH$. On a alors $\dot x^{-1} \dot\sigma \dot x = -\dot\sigma$. Il est clair que $\dot x \tilde \bM \dot x^{-1}$ et $\tilde\bM$ ne sont pas conjugu\'es dans $\bH$. 
Au contraire, l'image $x$ de $\dot x$ dans $\PSO_{4r}$ centralise l'image $\sigma$ de $\dot\sigma$. Par cons\'equent, les images de $\tilde \bM$ et de $x \tilde \bM x^{-1}$, qui restent non conjugu\'es, appartiennent tous les deux \`a $\cM$.

\subsection{Les groupes $\Spin$ et $\hSpin$}

Si $\bG$ est l'un des groupes $\Spin_m$ ou $\hSpin_{4m}$, posons $\bar\bG = \SO_m$ ou $\PSO_{4m}$, respectivement. $\bar\bG$ est donc un quotient de $\bG$ par un sous-groupe central de cardinal $2$. D\'efinissons $\bar\bL$, $\bar\bT_1$, 
$\bar\bH$, $\bar\Sigma_\semi$ et $\bar\bM_x$ comme \'etant les images dans 
$\bar\bG$ de $\bL$, $\bT_1$, $\bH$, $\Sigma_\semi$ et $\bM_x$. Le th\'eor\`eme \'etant d\'ej\`a \'etabli pour $\SO_m$ et $\PSO_{4m}$, on sait que les $\bar\bM_x$ se r\'epartissent en une ou deux classes de conjugaison sous $\bar\bH$.

Ainsi que nous l'avons remarqu\'e \`a la Section~\ref{subsect:m-graphe}, les questions de conjugaison de sous-groupes de L\'evi se r\'esolvent au niveau du graphe de Dynkin: deux sous-groupes de L\'evi sont conjugu\'es si et seulement si leur sous-graphes correspondants sont conjugu\'es par le groupe de Weyl. Les $\bM_x$ et les $\bar\bM_x$ \'etant des sous-groupes de L\'evi correspondants de $\bH$ et $\bar\bH$, on voit que la r\'epartition des $\bM_x$ en classes de conjugaison sous $\bH$ co\"incide avec celle des $\bar\bM_x$ sous $\bar\bH$.

\subsection{Le cas $\bG = \EE_6$, $\bL = \bM = (\AA_2)^2$}

Imposons les hypoth\`eses et prenons les notations de la
Remarque~\ref{rmk:m-conj-hyp}: on a $\sigma' = x^{-1}\sigma x = z\sigma$ pour un certain $z \in \bT_1$.

Le graphe de Dynkin compl\'et\'e $\tilde\Delta_\bG$ contient trois sous-graphes de type $(\AA_2)^2$:
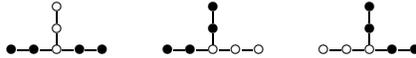
\[
\vcenter{\hbox{\small $\xymatrix@=4pt{
*{}&*{}& *{\circ} \ar@{-}[d] \\
*{}&*{}& *{\circ} \ar@{-}[d] \\
*{\bullet}\ar@{-}[r] & *{\bullet}\ar@{-}[r] &
*{\circ}\ar@{-}[r] & *{\bullet}\ar@{-}[r] & *{\bullet}}$}}
\qquad
\vcenter{\hbox{\small $\xymatrix@=4pt{
*{}&*{}& *{\bullet} \ar@{-}[d] \\
*{}&*{}& *{\bullet} \ar@{-}[d] \\
*{\bullet}\ar@{-}[r] & *{\bullet}\ar@{-}[r] &
*{\circ}\ar@{-}[r] & *{\circ}\ar@{-}[r] & *{\circ}}$}}
\qquad
\vcenter{\hbox{\small $\xymatrix@=4pt{
*{}&*{}& *{\bullet} \ar@{-}[d] \\
*{}&*{}& *{\bullet} \ar@{-}[d] \\
*{\circ}\ar@{-}[r] & *{\circ}\ar@{-}[r] &
*{\circ}\ar@{-}[r] & *{\bullet}\ar@{-}[r] & *{\bullet}}$}}
\]
Si $\Delta_\bH$ ne contient qu'un de ces trois graphes, alors le r\'esultat se d\'eduit des arguments de la Section~\ref{subsect:m-graphe}. Le seul sous-graphe 
propre de $\tilde\Delta_\bG$ qui en contient au moins deux est de type 
$(\AA_2)^3$, et celui-ci les contient tous les trois. Supposons d\'esormais que 
$\bH$ est de type $(\AA_2)^3$.

Le groupe $\bH$ est donc un quotient central de $(\SL_3)^3$. Notons $\mu_3$ le groupe des racines troisi\`emes de l'unit\'e, identifi\'e avec le groupe des 
matrices scalaires dans (\ie le centre de) $\SL_3$. Nous montrons maintenant que
\[
\bH \simeq (\SL_3 \times \SL_3 \times \SL_3)/\mu_3^\Delta,
\]
o\`u $\mu_3^\Delta$ est l'image du plongement diagonal 
$\mu_3 \hookrightarrow \mu_3 \times \mu_3 \times \mu_3$. Soit $K$ le noyau 
de l'application $(\SL_3)^3 \to \bH$. D'une part, on sait que $|\Cent(\bG)| = 3$,
et donc $|\Cent(\bL)/\Cento(\bL)|$ doit diviser $3$, mais si $K$ \'etait trivial 
(et donc $\bH \simeq (\SL_3)^3$), on pourrait en d\'eduire que 
$|\Cent(\bL)/\Cento(\bL)| = |\Cent(\SL_3 \times \SL_3)| = 9$. D'autre part, si 
$|K|$ valait $9$ ou $27$, il est facile de voir que $\bL$ serait \`a centre 
connexe, mais pour qu'un groupe de type $(\AA_2)^2$ admette un faisceau caract\`ere cuspidal, il ne doit pas \^etre \`a centre connexe. On conclut que $|K| = 3$, et sans perte de g\'en\'eralit\'e, on peut identifier $K$ avec $\mu_3^\Delta$.

Par un l\'eger abus de notation, nous allons \'ecrire des triplets $(a,b,c)
\in (\SL_3)^3$ pour d\'esigner des \'el\'ements de $\bH$. Soit $\omega \in
\mu_3$ une racine primitive troisi\`eme de l'unit\'e. Le centre de $\bH$
est l'ensemble $\{(a,b,c): a,b,c \in \mu_3\}$, o\`u on a, bien s\^ur,
l'identification $(\omega, \omega, \omega) = (1,1,1)$. 

Notons $\bT$ le groupe des matrices diagonales dans $\SL_3$, et soit $\bM_1$, 
$\bM_2$, et $\bM_3$ les images dans $\bH$ de $\SL_3 \times \SL_3 \times \bT$, 
$\bT \times \SL_3 \times \SL_3$, et $\SL_3 \times \bT \times \SL_3$, 
respectivement. Ces derniers sont des sous-groupes de L\'evi de $\bH$ correspondant aux trois graphes ci-dessus. Tout sous-groupe de L\'evi de $\bH$ de type $(\AA_2)^2$ est conjugu\'e \`a l'un de ces trois. Supposons, sans perte de g\'en\'eralit\'e, que $\bM = \bM_1$ et que $\bM_x$ est \'egal \`a l'un de $\bM_1$, $\bM_2$, $\bM_3$. Notre but est donc de d\'emontrer que $\bM_x = \bM_1$.

Soit $U$ une repr\'esentation irr\'eductible de $\bG$ de dimension $27$, et soit $V$ la repr\'esentation naturelle de $\SL_3$ de dimension $3$.  Alors
\[
U|_\bH \simeq 
\underbrace{V \otimes V^* \otimes \triv}_{U_1} \oplus
\underbrace{\triv \otimes V \otimes V^*}_{U_2} \oplus 
\underbrace{V^* \otimes \triv \otimes V}_{U_3}.
\]
On a $U_i = U^{\Cento(\bM_i)}$ pour $i = 1,2,3$. (Ici $U^{\Cento(\bM_i)}$ d\'esigne le sous-espace de $U$ sur lequel $\Cento(\bM_i)$ agit trivialement).

Un \'el\'ement $(a,b,c) \in \Cent(\bH)$ agit sur $U_1$ (resp.~$U_2$, $U_3$) par le scalaire $ab^{-1}$ (resp.~$bc^{-1}$, $ca^{-1}$). Nous pouvons maintenant identifier $\Cent(\bG)$ comme sous-groupe de $\Cent(\bH)$: c'est l'ensemble des \'el\'ements qui agissent sur $U$ par un scalaire,\ie l'ensemble des triplets $(a,b,c)$ o\`u $ab^{-1} = bc^{-1} = ca^{-1}$:
\[
\Cent(\bG) = \{(1,1,1),\ (1,\omega, \omega^2),\ (1,\omega^2,\omega)\}.
\]
Bien s\^ur, on a $\sigma \in \Cent(\bH)$; par contre, $\sigma \notin \Cent(\bG)$ (car $\bH = \Cento_\bG(\sigma)$). Il y a donc six possibilit\'es pour $\sigma$:
\begin{equation}\label{eqn:zh-liste}
\begin{aligned}
&(1,1,\omega),  &\qquad&(1,\omega,1),     &\qquad&(1,\omega^2,\omega^2) \\
&(1,1,\omega^2),&&(1,\omega,\omega),&&(1,\omega^2,1)
\end{aligned}
\end{equation}
Il est \`a noter que chaque \'el\'ement de $\Cent(\bH) \smallsetminus \Cent(\bG)$ agit sur les trois composantes $U_1$, $U_2$, $U_3$ par trois scalaires diff\'erents.

Montrons maintenant que $\sigma' \in \Cent(\bH) \smallsetminus \Cent(\bG)$ aussi. On sait que $\sigma' = z\sigma$ avec $z \in \bT_1$. Si $\sigma = (a,b,c)$, alors $\sigma' = (a,b,c')$ pour une certaine matrice $c' \in \bT$. \'Ecrivons
\[
\sigma' = \left(a,b,
\begin{bmatrix} c_1' &&\\ &c_2'&\\ &&c_3' \end{bmatrix}\right),
\qquad\text{o\`u $c_1'c_2'c_3' = 1$.}
\]
Cet \'el\'ement agit sur $U_1$ par le scalaire $ab^{-1}$, et ses valeurs propres sur $U_2$ (resp.~$U_3$) sont $bc_1'^{-1}$, $bc_2'^{-1}$, $bc_3'^{-1}$ (resp.~$c_1'a^{-1}$, $c_2'a^{-1}$, $c_3'a^{-1}$), chacune avec multiplicit\'e $3$. Mais ses valeurs propres doivent co\"incider avec celles de $\sigma$: on en d\'eduit imm\'ediatement que $c_1'$, $c_2'$, et $c_3'$ sont des racines troisi\`emes de l'unit\'e. Ensuite, la condition $c_1'c_2'c_3' = 1$ sur trois racines troisi\`emes de l'unit\'e implique qu'elles sont soit toutes \'egales, soit toutes distinctes. Mais si elles \'etaient toutes distinctes, $\sigma'$ aurait trois valeurs propres distinctes sur $U_2$, tandis que $\sigma$ n'en a que deux sur $U_2 \oplus U_3$. Il faut donc $c_1' = c_2' = c_3'$, \ie que $c'$ soit une matrice scalaire, et donc que $\sigma' \in \Cent(\bH)$.

Puisque $\sigma$ agit sur $U_1 = U^{\Cento(\bM)}$ par $ab^{-1}$, il faut que $\sigma'$ agisse sur $U^{\Cento(\bM_x)}$ par $ab^{-1}$. Mais $\sigma'$ agit sur $U_1$ par $ab^{-1}$, et par d'autres scalaires sur $U_2$ et $U_3$. On en d\'eduit que $U^{\Cento(\bM_x)} = U_1$, et donc que $\bM_x = \bM$.

\subsection{Le cas $\bG = \EE_7$, $\bL = \bM = (\AA_1)^3$}

Les hypoth\`eses et notations de la Remarque~\ref{rmk:m-conj-hyp} restent en vigueur.

Rappelons que $\EE_7$ contient plusieurs classes de conjugaison de sous-groupes de L\'evi de type $(\AA_1)^3$, dont une seule admet des faisceaux caract\`eres cuspidaux. Le graphe de Dynkin compl\'et\'e $\tilde\Delta_\bG$ contient deux sous-graphes correspondant \`a cette classe de conjugaison:
\[
\vcenter{\hbox{\small $\xymatrix@=4pt{
*{}&*{}&*{}& *{\bullet} \ar@{-}[d] \\
*{\bullet}\ar@{-}[r] & *{\circ}\ar@{-}[r] &
*{\bullet}\ar@{-}[r] & *{\circ}\ar@{-}[r] &
*{\circ}\ar@{-}[r] & *{\circ}\ar@{-}[r] & *{\circ}}$}}
\qquad
\vcenter{\hbox{\small $\xymatrix@=4pt{
*{}&*{}&*{}& *{\bullet} \ar@{-}[d] \\
*{\circ}\ar@{-}[r] & *{\circ}\ar@{-}[r] &
*{\circ}\ar@{-}[r] & *{\circ}\ar@{-}[r] &
*{\bullet}\ar@{-}[r] & *{\circ}\ar@{-}[r] & *{\bullet}}$}}
\]
Il s'ensuit que $\cM$ se r\'epartit en au plus deux classes de conjugaison sous $\bH$. Nous d\'emontrons par exemple maintenant que les membres de $\cM$ ne sont pas forc\'ement tous conjugu\'es.

Prenons pour $\bH$ l'unique sous-groupe (\`a conjugaison pr\`es) de type $\AA_3 \times \AA_1 \times \AA_3$. $\bH$ est donc un quotient de $\SL_4 \times \SL_2 \times \SL_4$. 

Soit $\bK_1$ l'image de $\SL_4 \times \SL_2 \times \bT$ dans $\bH$, et $\bK_2$ l'image de $\bT \times \SL_2 \times \SL_4$. Soit $\bM_1$ (resp.~$\bM_2$) le sous-groupe de L\'evi (unique \`a conjugaison pr\`es) de $\bK_1$ (resp.~$\bK_2$) de type $(\AA_1)^3$. Supposons, sans perte de g\'en\'eralit\'e, que $\bM = \bM_1$. Il est clair que $\bM_2$ n'est pas conjugu\'e \`a $\bM$ dans $\bH$.

Puisque $|\Cent(\bG)| = 2$, on sait que $|\Cent(\bK_1)/\Cento(\bK_1)| \le 2$. D'autre part, le fait que son sous-groupe de L\'evi $\bM_1$ admette un faisceau caract\`ere cuspidal implique que $|\Cent(\bK_1)/\Cento(\bK_1)| = 2$.  

Explicitons le centre de $\bH$. C'est un quotient de $\mu_4 \times \mu_2 \times \mu_4$. Des consid\'erations semblables \`a celles de la section pr\'ec\'edente permettent de trouver explicitement le noyau de cette application, en utilisant le fait que $|\Cent(\bK_1)/\Cento(\bK_1)| = 2$ et que le caract\`ere non trivial du centre de $\bK_1$ est \`a restriction non triviale sur chaque facteur quasi-simple de $\bM_1$. On trouve qu'on peut identifier
\[
\bH \simeq (\SL_4 \times \SL_2 \times \SL_4)/\mu_4^\Delta,
\]
o\`u $\mu_4^\Delta \subset \mu_4 \times \mu_2 \times \mu_4$ est le groupe cyclique engendr\'e par $(i,-1,i)$.

En particulier, on a $|\Cent(\bH)| = 8$. Deux \'el\'ements parmi les $8$ 
constituent $\Cent(\bG)$: \`a savoir, les triplets $(1,1,1)$ et $(-1,-1,-1)$. 
L'\'el\'ement $\sigma$ doit \^etre l'un des six \'el\'ements qui restent. 
\'Ecrivons un ensemble de repr\'esentants de ces six \'el\'ements:
\[
\text{ordre $4$:}\quad
\begin{array}{cc}
(1,1,i) & (1,-1,i) \\
(1,1,-i) & (1,-1,-i)
\end{array}
\qquad\qquad
\text{ordre $2$:}\quad
\begin{array}{c}
(1,1,-1) \\
(1,-1,1)
\end{array}
\]

Il est clair que $\bK_1$ et $\bK_2$ sont conjugu\'es sous $\bG$: leurs sous-graphes dans $\tilde\Delta_\bG$ sont conjugu\'es. Soit $x \in \bG$ un \'el\'ement tel que $x \bK_1 x^{-1} = \bK_2$ et $x \bK_2 x^{-1} = \bK_1$. La conjugaison par $x$ pr\'eserve $\bH$ et donc $\Cent(\bH)$. Puisqu'elle doit \'egalement pr\'eserver le facteur de type $\AA_1$ dans $\bH$, on voit que la conjugaison par $x$ doit stabiliser les deux \'el\'ements d'ordre $2$ dans $\Cent(\bH)$. (Il est \`a noter que $(1,1,-1) = (-1,1,1)$ dans $\Cent(\bH)$).

En r\'esum\'e, si l'on pose $\sigma = (1,1,-1)$ ou $\sigma = (1,-1,1)$, il existe un $x \in G$ qui stabilise $\sigma$, mais tel que $\bM_x = \bM_2$ n'est pas 
conjugu\'e dans $\bH$ \`a $\bM = \bM_1$.


\end{document}